\documentclass[12pt,twoside]{extarticle}

\usepackage[utf8]{inputenc}
\usepackage[english]{babel}
\usepackage{geometry}
\usepackage{fancyhdr}
\usepackage{titlesec}
\usepackage{amssymb}\usepackage{amsthm}\usepackage{amsmath}\usepackage{mathrsfs}\usepackage{mathtools}
\usepackage{yfonts}
\usepackage[]{enumitem}
\usepackage{comment}
\usepackage[dvipsnames]{xcolor}
\usepackage{microtype}
\usepackage{hyperref}
\usepackage{mathabx}
\newcommand{\claim}{\bigstar}
\renewcommand{\star}{\Asterisk}
\usepackage{mlmodern} 
\usepackage{moresize}
\usepackage{footmisc}

\setenumerate{itemsep=0pt,topsep=2pt,parsep=1pt}
\setitemize{itemsep = 0pt, topsep = 0pt, parsep = 0pt}
\hypersetup{colorlinks,linkcolor={red},citecolor={blue}}

\renewcommand{\thesection}{\arabic{section}}

\numberwithin{equation}{section}

\newcommand{\uc}{\MakeLowercase}
\newcommand{\nsiz}{\fontsize{12.5}{\baselineskip}\selectfont\ignorespaces}
\newcommand{\siz}{
\fontsize{13.5}{\baselineskip}\selectfont\ignorespaces
}
\newcommand{\sizz}{
\fontsize{15}{\baselineskip}\selectfont%
}

\titleformat{\chapter}
    {\Large \centering \sc \vspace*{0cm}}
    {Chapter \thechapter:}
    {5pt}
    {}
    [\vspace{1cm} \titlerule]
\titleformat{\section}
    {\large \centering \bfseries} 
    {\thesection.} 
    {5pt} 
    {} 
\titleformat{\subsection}
    [runin]
    {\bfseries}
    {\thesubsection}
    {5pt}
    {}
\renewenvironment{proof}{{\itshape Proof.}}{\hfill\qedsymbol\vspace{8pt}}

\newtheoremstyle{defbfnote}
    {8pt} 
    {8pt} 
    {} 
    {12pt} 
    {\scshape} 
    {.} 
    {.5 em} 
    {\thmname{\ignorespaces #1}\thmnumber{ #2}\thmnote{ \uc{(\ignorespaces\siz{#3})}}}
\newtheoremstyle{bfnote}
    {8pt} 
    {8pt} 
    {\itshape} 
    {12pt} 
    {\scshape} 
    {.} 
    {.5 em} 
    {\thmname{\ignorespaces #1}\thmnumber{ #2}\thmnote{ \uc{(\ignorespaces\siz{#3})}}}

\theoremstyle{bfnote}
\newtheorem{theorem}{\uc{\siz{Theorem}}}[section]
\newtheorem{corollary}[theorem]{\uc{\siz{Corollary}}}
\newtheorem{lemma}[theorem]{\uc{\siz Lemma}}
\newtheorem{proposition}[theorem]{\uc{\siz Proposition}}

\theoremstyle{defbfnote}
\newtheorem{definition}[theorem]{\uc{\siz{Definition}}}
\newtheorem{remark}[theorem]{\uc{\siz Remark}}

\setlist[enumerate,1]{label={(\roman*)}}
\setlist[itemize,1]{label = {$-$}}

\usepackage{amssymb,tikz}
\newcommand{\mysetminusD}{\hbox{\tikz{\draw[line width=0.6pt,line cap=round] (3pt,0) -- (0,8pt);}}}
\newcommand{\mysetminusT}{\mysetminusD}
\newcommand{\mysetminusS}{\hbox{\tikz{\draw[line width=0.45pt,line cap=round] (2pt,0) -- (0,5.5pt);}}}
\newcommand{\mysetminusSS}{\hbox{\tikz{\draw[line width=0.4pt,line cap=round] (1.5pt,0) -- (0,4pt);}}}

\renewcommand{\mathbb}{\mathbf}
\newcommand{\id}{\text{\normalfont id}}
\renewcommand{\P}{\mathcal{P}}
\newcommand{\B}{\mathcal{B}}
\newcommand{\R}{\mathbb{R}}

\newcommand{\Rn}{{{\mathbb{R}}^n}}
\newcommand{\C}{{{\mathscr{C}}}}
\newcommand{\Sp}{{{\mathscr{S}}}}
\renewcommand{\L}{{\mathscr{L}}}
\newcommand{\Ln}{{\mathscr{L}^n}}
\newcommand{\norm}[1]{\lVert#1\rVert}
\newcommand{\spt}{\text{\normalfont spt}}
\renewcommand{\d}{{\mathfrak{d}}}

\renewcommand{\div}{\text{\normalfont div}}

\renewcommand{\ae}{\text{a.e.}}
\renewcommand{\phi}{\varphi}

\newcommand{\Isop}{\text{\normalfont Isop}}
\newcommand{\tv}{\textsc{tv}}
\newcommand{\m}{{n-1}}

\renewcommand{\epsilon}{\varepsilon}

\newcommand{\ac}{\text{\normalfont AC}}

\newcommand{\Lip}{\text{\normalfont Lip}}
\newcommand{\s}{{\P^p_\infty(\Rn)}}

\newcommand{\dRn}{\text{\normalfont d}_\Rn}

\newcommand{\pl}{\mathrm{PL}}
\newcommand{\diam}{\text{\normalfont diam}}
\newcommand{\Graph}{\text{\normalfont Graph}}
\newcommand{\BF}{\mathbf}
\newcommand{\D}{\hat{\mathfrak{d}}_\infty^p}

\makeatletter
\def\thanks#1{\protected@xdef\@thanks{\@thanks
        \protect\footnotetext{#1}}}
\makeatother

\title{\uppercase{\bfseries \large Absolutely continuous curves in spaces of compactly supported densities }}
\author{{\sc pietro aldrigo\textsuperscript{\textdagger}} \thanks{\textsuperscript{\textdagger}{\sc Universit\"at Bern, Mathematisches Institut (MAI), Sidlerstrasse 12, 3012 Bern, Schweiz}. Email address: \texttt{pietro.aldrigo@unibe.ch}}
\thanks{The author is supported by the Swiss National Science Foundation, grant number 200021-228012. }
}
\date{ }

\begin{document}

\newcommand{\minus}{\mathbin{\mathchoice{\mysetminusD}{\mysetminusT}{\mysetminusS}{\mysetminusSS}}}
\renewcommand{\ge}{>}
\renewcommand{\le}{<}
\renewcommand{\div}{\text{\normalfont div}}
\newcommand{\Lcorner}{\mathrel{\makebox[7pt][c]{\rule{0.4pt}{6.75pt}\rule{5.5pt}{.4pt}}}}
\renewcommand{\llcorner}{{\Lcorner}}
\renewcommand{\Tilde}{\widetilde}
\renewcommand{\hat}{\widehat}
\newcommand{\f}{f}

\renewcommand{\star}{\Asterisk}
\newcommand{\curlybrace}[1]{\left\{#1\right\}}

\maketitle

\begin{abstract}
	We give a constructive proof for existence of absolutely continuous curves connecting each pair $\mu,\nu\in \pl_\infty^p(\Rn)$, for every $1\leq p\leq \infty$, where $(\pl_\infty^p(\Rn),\d_\infty^p)$ is the complete metric space of absolutely continuous measures with density in $L^p(\Rn)$ and bounded support introduced in \cite{AldrigoPLpq}. 
\end{abstract}

\section{Introduction}

Optimal transport gives a natural way of measuring distances between probability measures by taking into account how mass has to be moved; we refer to \cite{FigalliGlaudo2021, Santambrogio2015, Villani2003,Villani2009} for general introductions to the subject. This point of view has proved particularly effective in the study of variational problems and evolution equations, where Wasserstein distances interact naturally with variational and PDE techniques; see, for instance, \cite{ambrosio2005gradient,Brasco2010,Lisini2007}. One of the main reasons for this success is the dynamical interpretation of Wasserstein distances, beginning with the Benamou--Brenier formula \cite{BenamouBrenier2000}. In its classical form, this formula gives a dynamic characterization of the $2$-Wasserstein distance, and it extends to $(\P_q(\Rn),W_q)$ for every $1<q\leq \infty$; see, for example, Chapter 17 in \cite{ambrosio2021} and appendix of \cite{AldrigoPLpq} for the limit case $q=\infty$.

The Benamou--Brenier formula establishes the geodesic metric structure of $(\P_q(\Rn),W_q)$ and underlies the formal Riemannian viewpoint on Wasserstein spaces known as Otto calculus \cite{Otto2001}, where absolutely continuous curves of measures are described through the continuity equation. This has become a basic tool in the theory of gradient flows and minimizing movements in Wasserstein spaces: the general minimizing movement scheme goes back to De Giorgi \cite{DeGiorgi1993} and Ambrosio \cite{Ambrosio1995}, while its Wasserstein formulation was initiated by Jordan, Kinderlehrer and Otto \cite{JKO}; see also \cite{ambrosio2005gradient}.

In \cite{AldrigoPLpq}, the author introduced the spaces $\pl_q^p(\Rn)$ in order to combine $q$-transport control with  control of $L^p$-norms of the densities. More precisely, for $(p,q)\in [1,\infty]\times(1,\infty]$, one considers probability measures $\mu=f\Ln$ with $f\in L^p(\Rn)$ and with finite $q$-moment, or bounded support when $q=\infty$, endowed with the metric
\begin{align*}
\d_q^p(f\Ln,g\Ln)
:=
W_q(f\Ln,g\Ln)
+
\norm{f-g}_{L^p(\Rn)}.
\end{align*}

The original motivation for introducing these spaces came from the minimizing movement scheme for the isoperimetric ratio $\Isop:\pl_\infty^\infty(\Rn)\to [0,\infty]$,
\begin{align*}
\Isop(f\Ln):=\frac{\norm{\delta f}}{\norm{f}_{L^{n/(n-1)}}},
\end{align*}
where $\norm{\delta f}$ is the total variation of $f$ and $n\geq 2$. In \cite{AldrigoPLpq}, it was proved that, in the space $\pl_\infty^\infty(\Rn)$, the implicit scheme
\begin{align*}
\mu^{\tau}_{k+1}\in \operatorname*{argmin}_{\substack{\mu\in \pl_\infty^\infty(\Rn)}}\left\{\Isop(\mu)+ \frac{1}{2\tau}\left(\d_\infty^\infty(\mu,\mu^\tau_{k})\right)^2\right\}
\end{align*}
admits generalized minimizing movements which are absolutely continuous with respect to $\d_\infty^\infty$. The same paper also gave a characterization of $\d_q^p$-absolutely continuous curves through weak solutions of the continuity equation satisfying a Sobolev-type condition on the velocity field.

There remains, however, a basic geometric question which is not answered by such a characterization. The characterization tells us when a curve is absolutely continuous, but it does not say whether curves with prescribed endpoints actually exist. The metric $\d_q^p$ is stronger than the Wasserstein distance, since it controls both the $q$-transport cost and the $L^p$-difference of the densities. Hence one cannot infer AC path-connectedness from the Wasserstein geometry alone. In principle, the space $\pl_q^p(\Rn)$ could have too few admissible curves, or even contain pairs of measures which cannot be joined by a finite-length path.

The aim of this paper is to prove that
 when $q=\infty$ this obstruction does not occur. More precisely, we show the following statement.

\begin{theorem}\label{th_dn}
	For every $1\leq p\leq \infty$, each pair $\mu,\nu\in \pl_\infty^p(\Rn)$ can be connected by an absolutely continuous arc in $\pl_\infty^p(\Rn)$.
	
\end{theorem}

The proof is constructive: given an arbitrary compactly supported probability density, we explicitly deform it into a fixed normalized Lebesgue measure on a cube.

The main part of the construction is one-dimensional. Starting from a measure in $\pl_\infty^p(\R)$, we first smooth it by a controlled convolution-interpolation-type deformation (Definition \ref{def_interpolution}). We then translate and rescale the resulting density, perform a dilation-cutoff procedure (see Definition \ref{def_dilation-cutoff}), and finally interpolate it with the uniform probability measure on $[-1,1]$. Each of these operations is shown to be $\d_\infty^p$-Lipschitz. The estimates are kept separate for the $W_\infty$ component and for the $L^p$ component of the metric. This gives a path between any compactly supported one-dimensional density and the measure $\frac{1}{2}\L^1\llcorner[-1,1]$. Reversing and concatenating two such connections gives a Lipschitz continuous (hence absolutely continuous) curve between arbitrary endpoints in $\pl_\infty^p(\R)$.

The higher-dimensional argument is obtained by induction. We write $\Rn=\R^{n-1}\times\R$ and disintegrate a measure with respect to its projection onto the first $n-1$ variables. The one-dimensional construction is then applied along the last-coordinate fibers. The uniform estimates obtained in dimension one allow us to transform the fibers while preserving the required $W_\infty$ and $L^p$ bounds. In this way, the measure is connected to the product of its marginal on $\R^{n-1}$ with the uniform probability measure on $[-1,1]$. The induction hypothesis is then applied to the marginal, leaving the last coordinate unchanged. This argument reduces the original measure to the normalized Lebesgue measure on the unit cube $[-1,1]^n$.

The paper is organized as follows. In Section 2 we recall the notation and the preliminary facts on absolutely continuous curves in metric spaces, the $W_\infty$ distance, the continuity equation, and the spaces $\pl_\infty^p(\Rn)$. Section 3 contains the one-dimensional construction and proves connectivity by $\d_\infty^p$-absolutely continuous paths of $\pl_\infty^p(\R)$. In Section 4 we prove the higher-dimensional result by induction.

\section{Preliminaries}

\subsection{Curves in metric spaces.}
{
\renewcommand{\d}{{\normalfont\texttt{d}}}

\begin{definition}
	Let $(\Sp,\d)$ be a metric space. Fix a compact interval $I:=[a,b]$  and let $1\leq r\leq \infty$. A curve $\omega:I\to \Sp$ is \emph{$(\d,r)$-absolutely continuous} (or simply \emph{$\d$-absolutely continuous}, if $r=1$, and $\d$\emph{-Lipschitz}, if $r=\infty$), and we will write $\omega\in \ac^r(I;\Sp)$ (resp. $\ac(I;\Sp)$, $\Lip(I;\Sp)$), if there exists a function $m\in L^r(I)$ such that
	\begin{align}\label{eq_def_AC_metric}
		\d(\omega(t),\omega(s))\leq \int_s^t m(r)\,dr\quad \forall s,t\in I,\,\,  s<t.
	\end{align}
\end{definition}

\begin{proposition}
	Let $(\Sp,\d)$ be a metric space, fix $1\leq r\leq \infty$ and let $\omega\in \ac^r([a,b];\Sp)$.
	\begin{enumerate}
	\item We have
	\begin{align}\label{def_metricDerivative}
	\exists \lim_{{h\to 0}}\frac{\d(\omega(t+h),\omega(t))}{|h|}=:|\dot\omega|(t)\quad \text{for almost every }a\leq t\leq b.
	\end{align}
	\item The function $|\dot\omega|$ defined almost everywhere by \eqref{def_metricDerivative} belongs to $L^r([a,b])$ and satisfies \eqref{eq_def_AC_metric}.
	\item  If $m\in L^r([a,b])$  satisfies \eqref{eq_def_AC_metric}, then $|\dot\omega|\leq m$.
	\end{enumerate}
\end{proposition}

The function $|\dot \omega|$ is called \emph{$\d$-metric derivative of $\omega$}, or simply \emph{metric derivative of $\omega$}, when there is no ambiguity on the underlying metric structure.

Let $(\mathscr{S},\d)$ be a metric space and $\gamma_j\in \ac^p([s_j,t_j];\mathscr{S})$ for every $j=1,\dots,N$ be curves with the property $\gamma_j(t_j) = \gamma_{j+1}(s_{j+1})$ for every $1\leq j\leq N-1$. Let $\tau_j:= t_j-s_j$ and define $T_j:=\sum_{m=1}^j\tau_m$ for every $j=1,\dots,N$. The concatenation of $\gamma_1,\dots,\gamma_N$ is the path $\gamma_1\diamond\cdots\diamond \gamma_N\in \ac^p([0,T_N];\mathscr{S})$ defined by 
\begin{align*}
	\gamma_1\diamond\cdots\diamond \gamma_N(t):= \gamma_j(t-T_{j-1}+s_j),\quad \text{ if }T_{j-1}\leq t\leq T_j\quad \forall 0\leq t\leq T_N.
\end{align*} 
}

\subsection{The $\infty$-Wasserstein space.}

Let $n\geq 1$ be a fixed dimension. We denote by $(\Rn,\dRn)$ the Euclidean $n$-dimensional metric space. Consider the projections $\pi_j:\Rn\times\Rn\to \Rn$ defined by $\pi_j(x_1,x_2):=x_j$ for $j\in \{1,2\}$. Given two probability measures $\mu,\nu\in \P(\Rn)$ we define the set of \emph{couplings of $\mu$ and $\nu$} as
\begin{align*}
	\Gamma(\mu,\nu):=\left\{\gamma\in \P(\Rn\times \Rn): (\pi_1)_\#\gamma  = \mu,\, (\pi_2)_\#\gamma = \nu\right\}.
\end{align*} 

\begin{definition}
	The \emph{$\infty$-Wasserstein space} is the metric space $(\P_\infty(\Rn),W_\infty)$, where $\P_\infty(\Rn)$ is the family of Borel probability measures with bounded support and
	\begin{gather*}
	 W_\infty(\mu,\nu):=\inf_{\gamma\in \Gamma(\mu,\nu)}\norm{\dRn(\cdot,\cdot)}_{L^\infty(\gamma)} .
	\end{gather*}
  
\end{definition}

\begin{definition}\label{def_WeakSol}
	A couple $(\mu_{(\cdot)},v_{(\cdot)})$ is called \emph{weak solution of the continuity equation in $[0,T]$} (or simply \emph{solution of the continuity equation}) if $\mu_{(\cdot)}:[0,T]\to \P(\Rn)$, $[0,T]\times \Rn \ni (t,x)\mapsto v_t(x)\in \Rn$ is a Borel map, the integrability condition
\begin{align}\label{def_contEqIntegrability}
\int_0^T\norm{v_t}_{L^1(\mu_t)}\,dt \le \infty
\end{align}	
holds
	 and the partial differential equation
	 \begin{align*}
	 \partial_t \mu_t + \div(v_t\mu_t) = 0
	 \end{align*}
	 is satisfied in the distributional sense, namely
	\begin{align}\label{def_contEq}
	\int_0^T\int_\Rn\,\left(\partial_t \phi(t,x) + \langle \nabla \phi(t,x),v_t(x)\rangle\right)\,d\mu_t(x)\,dt = 0 \quad \forall \phi\in \C^\infty_c((0,T)\times \Rn). 
	\end{align}
	In this case, we call $v_{(\cdot)}$ the \emph{velocity field of $\mu_{(\cdot)}$}.
\end{definition}

The following lemma corresponds to Theorem 4.5 in \cite{AldrigoPLpq} (see also Theorem 8.3.1 in \cite{ambrosio2005gradient} for the analogous statement in $P_q(\Rn)$ with $1<q<\infty$).

\begin{lemma}\label{lem_ineqA}
	Let $\mu_{(\cdot)}:[0,T]\to \P_\infty(\Rn)$ be a narrowly continuous solution of the continuity equation for a Borel velocity-field $(t,x)\mapsto v_t(x)$ such that the function $t\mapsto \norm{v_t}_{L^\infty(\mu_t)}$ belongs to $L^1([0,T])$. Then $\mu_{(\cdot)}\in \ac([0,T];\P_\infty(\Rn))$ and the $W_\infty$-metric derivative $t\mapsto |\dot\mu|^{W_\infty}(t)$ of $\mu_{(\cdot)}$ satisfies
	\begin{align*}
	|\dot\mu|^{W_\infty}(t)\leq \norm{v_t}_{L^\infty(\mu_t)}\quad \ae\,0\leq t\leq T.
	\end{align*}
\end{lemma}

Let $n=1$ and consider $\mu\in \P(\R)$. The \emph{cumulative distribution function of $\mu$} is the non-decreasing right-continuous function $F_\mu:\R\to [0,1]$ defined by
\begin{align*}
F_\mu(x):=\mu((-\infty,x])\quad \forall x\in \R.
\end{align*}
For any non-decreasing right-continuous function $F:\R\to [0,1]$ its \emph{pseudo-inverse function} is the map $F^{-1}:[0,1]\to \R$ defined as
\begin{align*}
F^{-1}(y):=\inf\left\{z\in \R: F(z)>y\right\}.
\end{align*} 

If $\mu,\nu\in \P(\R)$ and $\mu = f\L^1$ is absolutely continuous, the \emph{monotone rearrangement of $\mu$ into $\nu$} is the map Borel $T:= F_\nu^{-1}\circ F_\mu:\R\to\R$. It turns out (see e.g. Section 1.4 in \cite{FigalliGlaudo2021}) that $T$ is a transport map from $\mu$ to $\nu$ (i.e. $(\id\times T)_\#\mu \in \Gamma(\mu,\nu)$). Moreover, it is elementary to verify that for every $r,s>0$, the following implication holds
\begin{align}\label{eq_monotone}
F_\nu(x-r)\leq F_\mu(x)\leq F_\nu(x+s)\quad \forall x\in \R\implies -r\leq T(x)-x\leq s\quad \mu\text{-a.e. }x\in \R.
\end{align}

\subsection{The spaces $\pl_\infty^p(\Rn)$.}

\begin{definition}
	Let $1\leq p\leq \infty$. The space $(\pl_\infty^p(\Rn),\d_\infty^p)$ is the family
	\begin{align*}
	\pl_\infty^p(\Rn):=\left\{\mu=f\Ln\in \P_\infty(\Rn): \,f\in L^p(\Rn)\right\},
\end{align*}
	endowed	with the metric
	\begin{align*}
	\d_\infty^p(f\Ln,g\Ln):= W_\infty(f\Ln,g\Ln) + \norm{f-g}_{L^p}.
	\end{align*}
\end{definition}

If $f_{(\cdot)}\in \ac^r([a,b];L^p(\Rn))$, for some $1\leq r\leq \infty$, we denote by $|\dot f|^{L^p}:[a,b]\to [0,\infty]$ the $(L^p(\Rn),\norm{\cdot}_{L^p})$-metric derivative of $f_{(\cdot)}$. Similarly, we use the notation $|\dot \mu|^{W_\infty}$ for the $(\P_\infty(\Rn),W_\infty)$-metric derivative of a curve $\mu_{(\cdot)}\in \ac^r([a,b];\P_\infty(\Rn))$. If $\mu_{(\cdot)}=f_{(\cdot)}\Ln\in \ac^r([a,b];\pl_\infty^p(\Rn))$, then the $(\pl_\infty^p(\Rn),\d_\infty^p)$-metric derivative of $\mu_{(\cdot)}$ is simply denoted by $|\dot \mu|$, and it is easy to check that 
\begin{align*}
|\dot \mu|(t) = |\dot \mu|^{W_{\infty}}(t) + |\dot f|^{L^p}(t)
\end{align*}
for almost every $a\leq t\leq b$.

We say that a Borel function $F:[0,T]\times \Rn\to\R$ is \emph{$L^1$ in time and $L^p$ in space}, and write $F\in L^1_tL^p_x([0,T]\times \Rn)$ if for almost every $0\leq t\leq T$ the $t$-section $F_t:=F(t,\cdot)$ belongs to $L^p(\Rn)$ and
\begin{align*}
\int_0^T\norm{F_t}_{L^p}\,dt<\infty.
\end{align*}

\begin{theorem}\label{th_ACWqp_pNotOne}
	Let $1\leq p\leq \infty$. A curve $\mu_{(\cdot)}=f_{(\cdot)}\Ln:[0,T]\to \pl_\infty^p(\Rn)$ is $\d_\infty^p$-absolutely continuous if and only if $\mu_{(\cdot)}$ is $W_\infty$-absolutely continuous and $f_{(\cdot)}$ is $L^p$-absolutely continuous. In particular, if $\mu_{(\cdot)}$ is narrowly continuous and there exists a Borel time-dependent vector field $(t,x)\mapsto v_t(x)$ such that
	\begin{enumerate}
	\item \label{th_ACPinftyp_0} $\int_0^T\norm{v_t}_{L^1(\mu_t)}\,dt<\infty$,
	\item \label{th_ACPinftyp_i} the vector-valued function $(t,x)\mapsto v_t f_t(x)$ admits a weak $x$-divergence $\div(v f)\in L^1_tL^p_x([0,T]\times \Rn)$,
	\item \label{th_ACPinftyp_ii} the real-valued function $(t,x)\mapsto f_t(x)$ admits a weak $t$-derivative $\partial_t f\in L^1_tL^p_x([0,T]\times \Rn)$,
	\item \label{th_ACPinftyp_iibis} $\partial_t f + \div(vf) = 0$ in $L^1_tL^p_x([0,T]\times \Rn)$, and
	\item \label{th_ACPinftyp_iiii} the function $t\mapsto \norm{v_t}_{L^\infty(\mu_t)}$ belongs to $L^1([0,T])$,
	\end{enumerate}
	then $\mu_{(\cdot)}$ is $\d_\infty^p$-absolutely continuous and 
	\begin{align*}
	|\dot \mu|(t)\leq \norm{v_t}_{L^\infty(\mu_t)} + \norm{\div(v_tf_t)}_{L^p}\quad\text{for almost every }0\leq t\leq T.
	\end{align*} 
\end{theorem}

A proof of the above statement is found in \cite{AldrigoPLpq}.

\section{One-dimensional construction}
\newcommand{\partialf}{\mathcal{D}f}
\renewcommand{\d}{{\mathfrak{d}_\infty^p}}
\renewcommand{\L}{\mathscr{L}^1}

\renewcommand{\s}{{\pl_q^p(\Rn)}}

We introduce here the main definitions required for the construction of $\d$-Lipschitz paths that connect arbitrary measures in $\pl_\infty^p(\R)$. 

Recall that if $\mu,\nu\in \P(\R)$ are two arbitrary Borel probability measures, then the \emph{convolution of $\mu$ and $\nu$} is the Borel probability $\mu\star \nu\in \P(\R)$ defined by
\begin{align*}
\mu\star \nu(B):= \int_\R\nu(B-x)\,d\mu(x) = \int_\R \mu(B-y)\,d\nu(y)\quad \forall B\in \B(\R).
\end{align*}

\begin{definition}[interpolution]\label{def_interpolution}
 Let $\rho\in \C^\infty_c(\R)$ be a non-negative function with $\int_\R \rho = 1$ and let $\mu\in \P(\R)$ be a probability measure. The \emph{$\rho$-interpolution of $\mu$} is the curve $\mu_{(\cdot)}:[0,1]\to \P(\R)$ defined by $\mu_t := \mu \star \kappa_t$ for every $0\leq t\leq 1$, where
 \begin{align*}
 \kappa_t :=\begin{cases}
 	\delta_0&,\text{ if }t=0,\\
 	(1-t)\delta_0 + \rho(\cdot/t)\L&,\text{ if }0<t\leq 1
 \end{cases},
\end{align*}
with $\delta_0\in \P(\R)$ being the Dirac's delta concentrated on $0\in \R$.
 
\end{definition}

\begin{definition}[shift]\label{def_shift}
Let $w\in \R$ and $\mu\in \P(\R)$ be fixed. The \emph{$w$-shift of $\mu$} is the curve $\mu_{(\cdot)}:[0,1]\to \P(\R)$ defined by $\mu_t(B):=\mu(B-tw)$ for every $B\in \B(\R)$ and every $0\leq t\leq 1$.
\end{definition}

Whenever $B\in \B(\R)$ is a Borel subset and $\lambda \geq 0$ is a non-negative scalar, we define 
\begin{align*}
\lambda B :=\left\{\lambda b\in \R: b\in B\right\}.
\end{align*}

\begin{definition}[scaling]\label{def_scaling}
Let $\lambda \in (0,\infty)$ and $\mu\in \P(\R)$ be fixed. The \emph{$\lambda$-scaling of $\mu$} is the curve $\mu_{(\cdot)}:[0,1]\to \P(\R)$ defined by $\mu_t(B):=\mu((1+t(\lambda-1))^{-1}B)$ for every $B\in \B(\R)$ and every $0\leq t\leq 1$.
\end{definition}

\begin{definition}[dilation-cutoff]\label{def_dilation-cutoff}
Let $T>1$ and $\mu\in \P(\R)$ be a probability measure such that $\mu([-1/T,1/T])>0$. The \emph{$T$-dilation-cutoff of $\mu$} is the curve $\mu_{(\cdot)}:[1,T]\to \P(\R)$ defined by 
\begin{align*}
\mu_t(B):=\frac{\mu(t^{-1}B\cap [-1/t,1/t])}{\mu([-1/t,1/t])},
\end{align*}
for every $B\in \B(\R)$ and every $1\leq t\leq T$.
\end{definition}

\begin{definition}[flattening-interpolation]\label{def_flattening-interpolation}
Let $B\in \B(\R)$ be a Borel subset with positive and finite Lebesgue measure and let $\mu\in \P(\R)$. The \emph{$B$-flattening-interpolation of $\mu$} is the curve $\mu_{(\cdot)}:[0,1]\to \P(\R)$ defined by
\begin{align*}
	\mu_t := (1-t)\mu + t\frac{\L\llcorner B}{\L(B)}
\end{align*} 
for every $0\leq t\leq 1$.
\end{definition}

For every $\ell>0$, we denote by $\mathfrak{C}^\infty_\ell(\R)$ the family of functions $\rho\in \C^\infty_c(\R)$ that are even, with $\int_\R \rho = 1$, $\rho>0$ in $[0,\ell)$, $\rho'<0$ in $(0,\ell)$, and $\rho\equiv 0$ in $[\ell,\infty)$.

\begin{definition}[standard connection]\label{def_standardConnection}
Let $[a^*,b^*]$ be a compact interval and let $\mu\in \P_\infty(\R)$ be a probability measure such that $[\alpha,\beta]:=[\min\,\spt\,\mu,\max\,\spt\,\mu]\subseteq [a^*,b^*]$. Define  $r:=(b^*-a^*)/2$, $M:=\max\{|a^*|,|b^*|\}$, $a:=\alpha-3r$, $b:=\beta+3r$, $w:=-(a+b)/2$, $\lambda:=b-a$ and let $\rho\in \mathfrak{C}^\infty_{3r}(\R)$. 
The \emph{$(a^*,b^*,\rho)$-standard connection of $\mu$ with $[-1,1]$} is the curve $\mu_{(\cdot)}:[0,7]\to \pl_\infty^p(\R)$ obtained by the concatenation  $\mu^1_{(\cdot)}\diamond\cdots\diamond\mu^5_{(\cdot)}$, where $\mu^1_{(\cdot)}$ is the $\rho$-interpolution of $\mu$, $\mu^2_{(\cdot)}$ is the $w$-shift of $\mu^1_1$, $\mu^3_{(\cdot)}$ is the $(2/\lambda)$-scaling of $\mu^2_1$, $\mu^4_{(\cdot)}$ is the $4$-dilation-cutoff of $\mu^3_1$ and $\mu^5_{(\cdot)}$ is the $[-1,1]$-flattening-interpolation of $\mu^4_4$.
\end{definition}

Throughout this work, we adopt the following notation. Whenever $c>0$ is a positive number, we define
\begin{align*}
c^\frac{1}{\infty}:=1\quad \text{and}\quad c^\frac{\infty-1}{\infty}:=c.
\end{align*}
Moreover, for every $r>0$, we denote by $\rho^r\in \mathfrak{C}^\infty_{3r}(\R)$ the function 
\begin{align}\label{rho}
\rho^r(x):=\frac{e^{\frac{1}{(x^2-9r^2)}}}{\int_{-3r}^{3r}e^{\frac{1}{(y^2-9r^2)}}\,dy}\chi_{[-3r,3r]}(x).
\end{align}

\begin{theorem}\label{th_d1}
Fix $1\leq p\leq \infty$. Let $[a^*,b^*]\subseteq \R$ be a compact interval and set $r:=(b^*-a^*)/2$. Suppose $r>1/3$. Then there exist constants $L,M,m >0$ that depend only on $a^*$ and $b^*$ that satisfy the following property. For any $Q>0$, any probability measure $\mu\in \pl_\infty^p(\R)$ with $\spt\,\mu\subseteq [a^*,b^*]$ and $\norm{f}_{L^p}\leq Q$, the $(a^*,b^*,\rho^r)$-standard connection of $\mu$ with $[-1,1]$ is a $\d$-Lipschitz curve $\mu_{(\cdot)}=f_{(\cdot)}\L:[0,7]\to \pl_\infty^p(\R)$  and we have the estimates 
\begin{gather}\label{th_d1_1}
|\dot f|^{L^p}(t)\leq C^{L^p}_{a^*,b^*,Q}(t):=
\begin{cases}
	 2Q																&, \text{ if } 0\leq t < 1\\
	 LM(8r)^\frac{1}{p}												&, \text{ if } 1\leq t < 2\\
	 32 Lr^2(4r-1)(8r)^\frac{1}{p}									&, \text{ if } 2\leq t < 3\\
	 256\left(\frac{r^2 L}{m}\right)^2 2^\frac{1}{p}				 		&, \text{ if } 3\leq t < 6\\
	 \left(\frac{1}{2}+\frac{16r^2L}{m}\right)2^\frac{1}{p}			&, \text{ if } 6\leq t \leq 7
\end{cases},
\end{gather}
\begin{gather}
\label{th_d1_2}
|\dot \mu|^{W_\infty}(t)\leq C^{W_\infty}_{a^*,b^*}(t):=
\begin{cases}
	 3r 																		&, \text{ if } 0\leq t < 1\\
	 M																		&, \text{ if } 1\leq t < 2\\
	4r-1																		&, \text{ if } 2\leq t < 3\\
	2+\frac{16r^2 L}{m} 														&, \text{ if } 3\leq t < 6\\
	\left(\frac{1}{2} + \frac{16 r^2 L}{m}\right)\max\left\{\frac{2}{m},4\right\}	&, \text{ if } 6\leq t \leq 7
\end{cases}
\end{gather}
for almost every $0\leq t \leq 7$.
\end{theorem}

\begin{remark}
	Notice that the function $C^{W_\infty}_{a^*,b^*}$ defined in \eqref{th_d1_2} does not depend on the upper bound $Q$ of the $L^p$ norm on the density, nor on the value of $p$, but only on $a^*$ and $b^*$, while the  $C^{L^p}_{a^*,b^*,Q}$ defined in \eqref{th_d1_1} depends on the bound $Q$ only for $0\leq t< 1$.
\end{remark}

\begin{corollary}\label{cor}
	Let $1\leq p\leq \infty$. For every pair $\mu,\nu\in \pl_\infty^p(\R)$ there exists a $\d$-absolutely continuous path $\gamma_{(\cdot)}:[0,1]\to \pl_\infty^p(\R)$ such that $\gamma_0=\mu$ and $\gamma_1=\nu$.
\end{corollary}
\begin{proof}
By virtue of Theorem \ref{th_d1}, there are $\d$-Lipschitz paths $\mu_{(\cdot)}:[0,7]\to \pl_\infty^p(\R)$ and $\nu_{(\cdot)}:[0,7]\to \pl_\infty^p(\R)$ such that $\mu_0=\mu$, $\nu_0=\nu$ and $\mu_7=\nu_7=\frac{1}{2}\L\llcorner [-1,1]$. Then, define $\hat \nu_{(\cdot)}:[0,7]\to \pl_\infty^p(\R)$ as $\hat \nu_t := \nu_{7-t}$. The path $\gamma_{(\cdot)}:[0,1]\to \pl_\infty^p(\R)$ defined by $\gamma_t:= (\mu_{(\cdot)}\diamond \hat \nu_{(\cdot)})_{14t}$ is a $\d$-Lipschitz (hence a $\d$-absolutely continuous) curve and we have $\gamma_0= \mu$ and $\gamma_1 = \nu$.

\end{proof}

For the rest of this section, we tacitly assume that $p$ is a fixed exponent in $[1,\infty]$.

\begin{lemma}[lipschitz continuity of the interpolution]\label{lem_interp}
	Let $\mu\in \pl_\infty^p(\R)$ and denote by $[\alpha,\beta]:=[\min\,\spt\,\mu,\max\,\spt\, \mu]$. Let $r>(\beta-\alpha)/2$ and fix a function $\rho\in \mathfrak{C}^\infty_{3r}(\R)$. Let $\mu_{(\cdot)}=f_{(\cdot)}\L:[0,1]\to \pl_\infty^p(\R)$ be the $\rho$-interpolution of $\mu$. Then $\mu_{(\cdot)}$ is $\d$-Lipschitz and in particular
	\begin{align}\label{lem_interp_1}
	|\dot f|^{L^p}(t)\leq 2\norm{f}_{L^p}\quad\text{and}\quad |\dot\mu|^{W_\infty}(t)\leq 3r
	\end{align}
	for almost every $0\leq t\leq 1$. Moreover, the function $f_1$ satisfies:
	\begin{enumerate}
	\item \label{lem_interp_i} $f_1\in \C^\infty_c(\R)$, with $\spt\,f_1= [\alpha -3r,\beta + 3r]$ and $f_1>0$ in $(\alpha-3r,\beta+3r)$;
	\item \label{lem_interp_ii} $\Lip(f_1)\leq \Lip(\rho)$;
	\item \label{lem_interp_iii} $f_1$ is increasing in $(\alpha-3r,\alpha)$ and decreasing in $(\beta,\beta+3r)$.
	\end{enumerate}
\end{lemma}

\begin{proof}
Let $\mu$, $\alpha$, $\beta$, $r$ and $\rho$ as in the statement and define $\rho_t=\rho(\cdot/t)t^{-1}$. Define the convolution densities $g_t:= f\star \rho_t$ for every $0<t\leq 1$ and the measures $\nu_t:=g_t\L$, so that the $\rho$-interpolution of $\mu$ writes
\begin{align*}
\mu_t = \begin{cases} \mu &,\text{ if }t=0\\ (1-t)\mu + t\nu_t &,\text{ if }0<t\leq 1\end{cases} = \begin{cases}f\L&,\text{ if }t=0\\ ((1-t)f+tg_t)\L&,\text{ if }0<t\leq 1\end{cases}.
\end{align*}

First we show that the curve of densities $f_{(\cdot)}$ is $L^p$-Lipschitz continuous. Since $\norm{g_t}_{L^p}\leq \norm{f}_{L^p}$ for every $1\leq p\leq \infty$, we deduce that 
\begin{align*}
\norm{f_t-f_0}_{L^p} = t\norm{g_t-f}_{L^p}\leq 2t\norm{f}_{L^p}\quad \forall 0< t\leq 1.
\end{align*}
Moreover,
 the function $t\mapsto f_t(x)$ is differentiable in $(0,1]$ for almost every $x\in \R$ and
\begin{align*}
\partial_tf_t(x) = -f(x) + f\star \zeta_t(x),
\end{align*}
where $x\mapsto \zeta(x):=-\rho'(x)x$ is a smooth, even, non-negative and compactly supported function and $\zeta_t(x):=\zeta(\cdot/t)t^{-1}$. From $\norm{\zeta_t}_{L^1} = \int_\R \rho(y)\,dy =1$  it follows that
\begin{align*}
\norm{\partial_t f_t}_{L^p}\leq \norm{f}_{L^p} + \norm{f\star \zeta_t}_{L^p}\leq  2\norm{f}_{L^p}.
\end{align*}
 Whence, writing $f_t(x)-f_s(x)= \int_s^t\partial_tf_r(x)\,dr$ for almost every $x\in \R$ and using the Minkowski inequality for $L^p$-norms we obtain
\begin{align*}
\norm{f_t-f_s}_{L^p}\leq 2(t-s)\norm{f}_{L^p}\quad \forall 0<s\leq t\leq 1,
\end{align*} 
and thus the $L^p$-Lipschitz continuity of $t\mapsto f_t$ is proved.

We shall now show that for every given $0\leq t <t+h\leq 1$, if $T:\R\to\R$ is the monotone rearrangement of $\mu_t$ into $\mu_{t+h}$, then 
\begin{align}\label{interpolution_0}
|T(x)-x|\leq 3rh \quad \text{for }\mu_t\text{-almost every } x\in \R.
\end{align}
This in turn proves the $W_\infty$-Lipschitz continuity of $t\mapsto \mu_{t}$, and more precisely that
\begin{align}\label{eq_C2}
W_\infty(\mu_{t+h},\mu_t)\leq 3rh.
\end{align}
Define $r^*:=3r$ and the probability measures
\begin{align*}
 \kappa_t :=\begin{cases}
 	\delta_0&,\text{ if }t=0,\\
 	(1-t)\delta_0 + \rho(\cdot/t)\L&,\text{ if }0<t\leq 1
 \end{cases},
\end{align*}
with $\delta_0\in \P(\R)$ being the Dirac's delta concentrated on $0\in \R$,
and consider the cumulative distribution functions
\begin{gather}
	R(x):=\int_{-\infty}^x\rho(y)\,dy,\notag
	\quad R_t(x):=\int_{-\infty}^x\rho_t(y)\,dy = R(x/t),
	\\ K_t(x):=\kappa_t((-\infty,x]) = (1-t)\chi_{[0,\infty)}(x)+tR(x/t),\notag
	\\ F_t(x):=\mu_t((-\infty,x]) = \int_\R K_t(x-y)\,d\mu(y). \label{eq_Ft}
\end{gather}
Then, by definition of monotone rearrangement, we have that $T=F_{t+h}^{-1}\circ F_t$, where $F_{t+h}^{-1}$ is the pseudo-inverse of $F_{t+h}$.

We claim that 
\begin{align}\label{interpolution_1}
K_{t+h}(x-r^*h)\leq K_t(x) \leq K_{t+h}(x+r^*h)\quad \forall x\in \R.
\end{align}
Indeed, recalling that $\spt\,\rho = [-r^*,r^*]$ we deduce
\begin{align*}
K_{t+h}(x+r^*h) - K_t(x) \quad \begin{cases}
	\geq (t+h)R\left(\frac{x+r^*h}{t+h}\right)- tR\left(\frac{x}{t}\right)&,\text{ if }-\infty<x<0,\\
	= -h + (t+h)R\left(\frac{x+r^*h}{t+h}\right) - tR\left(\frac{x}{t}\right)&,\text{ if }0\leq x\leq r^*t\\
	= 0 &,\text{ if }r^*t<x<\infty
\end{cases}.
\end{align*}
Observe that
\begin{gather}
\frac{x+r^*h}{t+h}\geq \frac{x}{t}\quad \forall -\infty < x\leq r^*t,\label{pf_lem_iterp_1}
\\ 
R\text{ is non-decreasing in }\R\text{ and }\rho\text{ is non-increasing in }[0,\infty),\label{pf_lem_iterp_2}
\\
\frac{d}{dx}\left(-h + (t+h)R\left(\frac{x+r^*h}{t+h}\right) - tR\left(\frac{x}{t}\right)\right) = \rho\left(\frac{x+r^*h}{t+h}\right)-\rho\left(\frac{x}{t}\right).\label{pf_lem_iterp_3}
\end{gather} 
From \eqref{pf_lem_iterp_1} and \eqref{pf_lem_iterp_2} follows that
\begin{align*}
 K_{t+h}(x+r^*h)-K_t(x) \geq 0 \quad \forall -\infty < x <0.
\end{align*}
Finally, from \eqref{pf_lem_iterp_3} and  absolute continuity, we deduce that $x\mapsto K_{t+h}(x+r^*h)-K_t(x)$ is non-increasing in $[0,r^*t]$ and vanishing at $x=r^*t$. Therefore, the second inequality of \eqref{interpolution_1} is proved. The first inequality of \eqref{interpolution_1} is proved analogously.

As a consequence of \eqref{eq_Ft} and \eqref{interpolution_1}, we deduce that
\begin{gather}
F_{t+h}(x-r^*h)\leq F_t(x)\leq F_{t+h}(x+r^*h)\quad \forall x\in \R.\label{interpolution_3}
\end{gather}
Combining \eqref{eq_monotone} with \eqref{interpolution_3} yields \eqref{interpolution_0}.

We now prove properties \ref{lem_interp_i}, \ref{lem_interp_ii} and \ref{lem_interp_iii}. By classical properties of the convolution, the support of $f_1$ is the $3r$-neighborhood of $\spt\,f$. Since $2r>\diam \,\spt\,\mu$, $\spt\,f_1 = [\alpha-3r,\beta+3r]$. Fix $x\in (\alpha-3r,\beta+3r)$, then $(x-3r,x+3r)\cap (\alpha,\beta)$ is a non-trivial open interval that contains at least one element of $\{\alpha,\beta\}$. Thus  
\begin{align*}
f_1(x) = \int_{x-3r}^{x+3r}f(z)\rho(x-z)\,dz >0.
\end{align*} 
Property \ref{lem_interp_ii} is a standard property of the convolution. Finally, if $\beta <x<y<\beta+3r$, then, recalling that $\spt\,f\subseteq [\alpha,\beta]$ and that $\rho$ is strictly decreasing in $(0,3r)$, we deduce
\begin{align*}
	f_1(x)= \int_{x-3r}^{\beta} f(z)\rho(x-z)\,dz > \int_{y-3r}^{\beta} f(z)\rho(y-z)\,dz = f_1(y).
\end{align*}
Similarly, one shows that $f_1$ is increasing in $(\alpha-3r,\alpha)$. 
This concludes the proof of \ref{lem_interp_iii}.

\end{proof}

\begin{lemma}[lipschitz continuity of the shift]\label{lem_shift}
Let $\mu=f\L\in \pl_\infty^p(\R)$ be such that $f\in \Lip(\R)$ and $\spt\,f\subseteq [a^*,b^*]$, and fix $w\in \R$. Then the $w$-shift $\mu_{(\cdot)}=f_{(\cdot)}\L$ of $\mu$ is $\d$-Lipschitz and, in particular, the estimates
\begin{align}\label{lem_Lip_1}
	|\dot f|^{L^p}(t)\leq \Lip(f)|w|(b^*-a^*)^\frac{1}{p}\quad \text{and}\quad |\dot\mu|^{W_\infty}(t)\leq |w|
\end{align}
hold for almost every $0\leq t\leq 1$.
\end{lemma}

\begin{proof}
	A (representative of the) density of the $w$-shift $\mu_t$ is $f_t(x) =f(x-tw)$ for every $0\leq t\leq 1$. Therefore $t\mapsto f_t(x)$ is Lipschitz for every $x\in \R$	and
	\begin{align*}
	|\partial_t f_t(x)| = |f'(x-tw)w| \leq \Lip(f)|w|\chi_{[a^*,b^*]}(x-tw)
	\end{align*}
	for almost every $(t,x)\in [0,1]\times \R$.
	Therefore, given $0\leq s<t\leq 1$, writing $f_t(x)-f_s(x)=\int_s^t\partial_t f_r(x)\,dr$ for almost every $x\in \R$ we deduce
	\begin{align*}
	\norm{f_t-f_s}_{L^p}\leq \int_s^t\norm{\partial_tf_r}_{L^p}\,dr\leq (t-s)\Lip(f)|w|(b^*-a^*)^\frac{1}{p}\quad \forall 0\leq t\leq 1.
	\end{align*}
	Hence, the first inequality of \eqref{lem_Lip_1} follows.
	
	Finally, observe that the time-dependent (1-dimensional) vector field $(t,x)\mapsto v_{t}(x)\equiv w$ defined on $[0,1]\times \R$ is a velocity field for $\mu_{(\cdot)}$. Therefore, Lemma \ref{lem_ineqA} implies the second inequality of \eqref{lem_Lip_1}.
	
\end{proof}

\begin{lemma}[lipschitz continuity of the scaling]\label{lem_scaling}
Let $\mu=f\L\in \pl_\infty^p(\R)$ be such that $f\in \Lip(\R)$ and $\spt\,f = [-\lambda/2,\lambda/2]$, with $\lambda \geq 2$. Then the $(2/\lambda)$-scaling $\mu_{(\cdot)}=f_{(\cdot)}\L$ of $\mu$ is $\d$-Lipschitz and, in particular, we have the estimates
\begin{align}\label{lem_scaling_1}
	|\dot f|^{L^p}(t)\leq \Lip(f)\frac{\lambda^2(\lambda-2)}{4}\lambda^\frac{1}{p}\quad \text{and}\quad |\dot \mu|^{W_\infty}(t)\leq \frac{\lambda}{2}-1 
\end{align}
for almost every $0\leq t\leq 1$.
\end{lemma}

\begin{proof}
	Let $\ell:= 2/\lambda\leq 1$. A (representative of the) density of the $\ell$-scaling at time $0\leq t\leq 1$ is 
	\begin{align*}
	f_t(x)= \frac{f\left(\frac{x}{1+t(\ell -1)}\right)}{1+t(\ell -1)}\quad \forall x\in \R.
	\end{align*}
	Since $1+t(\ell -1)\geq \ell$ for every $0\leq t\leq 1$, the function $(t,x)\mapsto f_t(x)$ is Lipschitz in $[0,1]\times \R$, and its $t$-derivative at almost every $(t,x)$ is
	\begin{align*}
	\partial_t f_t(x) = \frac{f'\left(\frac{x}{1+t(\ell-1)}\right)\frac{x}{1+t(\ell-1)} + f\left(\frac{x}{1+t(\ell-1)}\right)}{(1+t(\ell-1))^2}(1-\ell).
\end{align*}		
	It is easy to check that
	\begin{align*}
	|\partial_t f_t(x)|& 
	\leq \Lip(f)\frac{\lambda^2(\lambda-2)}{4}\chi_{[-\lambda/2,\lambda/2]}(x)\quad \text{a.e. }0\leq t\leq 1\quad \forall x\in \R.
	\end{align*}
	Thus, arguing as in the previous proofs, the first inequality of \eqref{lem_scaling_1} follows.
	
	Define now the time-dependent (1-dimensional) vector field $(t,x)\mapsto v_t(x)$ 
	\begin{align*}
	v_t(x):= \frac{\ell-1}{1+t(\ell-1)}x.
	\end{align*}
	Then $(t,x)\mapsto v_t(x)f_t(x)$ is Lipschitz in $[0,1]\times \R$ and $\frac{d}{dx}(v_tf_t)(x) = -\partial_t f_t(x)$. Therefore $(\mu_{(\cdot)},v_{(\cdot)})$ is a solution of the continuity equation and, by Lemma \ref{lem_ineqA}, and using the fact that $|x|/(1+t(\ell-1))\leq \lambda/2$ in $\spt\,\mu_t$, we obtain
	\begin{align*}
	|\dot\mu|^{W_\infty}(t) \leq \norm{v_t}_{L^\infty(\mu_t)}\leq (1-\ell)\frac{\lambda}{2}=\frac{\lambda}{2}-1.
	\end{align*}
\end{proof}

\begin{lemma}[lipschitz continuity of the dilation-cutoff]\label{lem_dilCut}
	Let $\mu=f\L\in \pl_\infty^p(\R)$ with $f\in \C^\infty(\R)$, $spt\,f =[-1,1]$, $f>0$ in $(-1,1)$. Suppose that $f$ is increasing in $(-1,-1/4)$, decreasing in $(1/4,1)$ and that there exists a constant $m>0$ such that $\min_{[-1/4,1/4]} f\geq 2m$. Let $\mu_{(\cdot)}=f_{(\cdot)}\L:[1,4]\to \pl_\infty^p(\R)$ be the $4$-dilation-cutoff of $\mu$. Then $\mu_{(\cdot)}$ is $\d$-Lipschitz and 
	\begin{align}\label{lem_dilCut_1}
	|\dot f|^{L^p}(t)\leq 4\frac{\norm{f}_{L^\infty}\Lip(f)}{m^2}2^\frac{1}{p}\quad \text{and}\quad |\dot \mu|^{W_\infty}(t)\leq 2+\frac{\norm{f}_{L^\infty}}{m}
	\end{align}
	hold for almost every $1\leq t\leq 4$.
\end{lemma}

\begin{proof}
{
	Let $\mu_{(\cdot)}$ be the $4$-dilation-cutoff of $\mu$. Then $\mu_{(\cdot)}=f_{(\cdot)}\L$, where 
	\begin{align*}
	f_t(x)=\frac{f(x/t)\chi_{[-1,1]}(x)}{\int_{-1}^1f(y/t)\,dy}\quad \forall (t,x)\in [1,4]\times \R.
\end{align*}	 
	Since $f$ is smooth in $\R$, then so it is $(t,x)\mapsto f_t(x)$ in $[1,4]\times (-1,1)$. Observe that
	\begin{align*}
	-\partial_tf_t(x) = \frac{f'(x/t)x\int_{-1}^1 f(y/t)\,dy - f(x/t)\int_{-1}^1f'(y/t)y\,dy}{\left(t\int_{-1}^1f(y/t)\,dy\right)^2}\quad \forall (t,x)\in [1,4]\times (-1,1).
	\end{align*}
	Define $(t,x)\mapsto v_t(x)$ in $[1,4]\times \R$ as
	\begin{align*}
	v_t(x):=-\frac{\int_{-1}^x \partial_t f_t(z)\,dz}{f_t(x)}\chi_{[-1,1]}(x).
	\end{align*}
	Since $f>0$ in $(-1,1)$, the function $f_t$ is strictly positive in $(-1,1)$ for every $1\leq t\leq 4$. Moreover $z\mapsto \partial_t f_t(z)$ is integrable in $(-1,1)$ and $\int_{-1}^1\partial_t f_t(z)\,dz=0$ for every $1\leq t\leq 4$. Therefore, $(t,x)\mapsto v_t(x)$ is well-defined in $[1,4]\times\R$ and it is Borel-measurable.
	
	Let $\psi\in \C^\infty_c((1,4)\times \R)$ be an arbitrary test function. Then,  observing that $\int_{-1}^1\partial_tf_t(z)\,dz=0$, we deduce
	\begin{align}\label{eq_previous}
	\begin{split}
	\int_1^4\int_\R &\left(\partial_t\psi(t,x) + \partial_x\psi(t,x)v_t(x)\right)\,d\mu_t(x)\,dt\\ 
	&= \int_{-1}^1\int_1^4 \partial_t\psi_t f_t(x)\,dt\,dx - \int_1^4\int_{-1}^1 \partial_x \psi(t,x)\left( \int_{-1}^x \partial_t f_t(z)\,dz\right)\,dx\,dt\\
	&= -\int_{-1}^1\int_1^4 \psi(t,x)\partial_t f_t(x)\,dt\,dx + \int_1^4 \int_{-1}^1\psi(t,x)\partial_tf_t(x)\,dx\,dt\\
	&=0.
	\end{split}
	\end{align}
	Thus, $(\mu_{(\cdot)},v_{(\cdot)})$ is a solution of the continuity equation in $[1,4]$ in the sense of Definition \ref{def_WeakSol}.
	
	Recalling $0<m\leq \frac{1}{2}\min_{[-1/4,1/4]}f$, we deduce
	\begin{align}\label{eqqq_1}
	\int_{-1}^1 f(y/t)\,dy = t\int_{-1/t}^{1/t}f(\xi)\,d\xi \geq \int_{-1/4}^{1/4}f(\xi)\,d\xi \geq m.
	\end{align}
	Therefore,
	\begin{align*}
	|\partial_t f_t(x)| =\frac{\left|f'(x/t)x\int_{-1}^1 f(y/t)\,dy - f(x/t)\int_{-1}^1f'(y/t)y\,dy\right|}{\left(t\int_{-1}^1f(y/t)\,dy\right)^2} \leq 4\frac{\norm{f}_{\C^0}\norm{f'}_{\C^0}}{m^2}
	\end{align*}
	for every $(t,x)\in [1,4]\times(-1,1)$. Hence, by the standard integration argument, the first inequality of \eqref{lem_dilCut_1} is proved
	and the curve $t\mapsto f_t$ is $L^p$-Lipschitz.
	
	To estimate the $L^\infty$-norm of $v_t$, let us begin with some simple arithmetics. For $-1<x<1$, integrating by parts we deduce
	\begin{align*}
	\begin{split}
	v_t(x) 
	&= \frac{\int_{-1}^x f'(z/t)z\,dz}{f(x/t)t^2} - \frac{\int_{-1}^x f(z/t)\,dz\, \int_{-1}^{1}f'(y/t)y\,dy}{f(x/t)t^2\int_{-1}^1f(y/t)\,dy}\\
	&=\frac{x}{t}+ \frac{(1-\omega_t(x))}{tf(x/t)}f(-1/t) -\frac{\omega_t(x)}{t f(x/t)}f(1/t).
	\end{split}
	\end{align*}
	where
	\begin{align*}
	\omega_t(x):=\frac{\int_{-1}^x f(z/t)\,dz}{\int_{-1}^1 f(y/t)\,dy}\in [0,1]\quad \forall -1\leq x\leq 1.
	\end{align*}
	Therefore, defining $A_{-1}^t,A_1^t:(-1,1)\to [0,\infty)$ as
	\begin{align*}
	A_{-1}^t(x):= \frac{(1-\omega_t(x))}{tf(x/t)}f(-1/t)\quad\text{and}\quad A_1^t(x):=\frac{\omega_t(x)}{t f(x/t)}f(1/t)
	\end{align*}
	for every $1\leq t\leq 4$,
	we infer the estimate
	\begin{align}\label{pf_lem_ext&cut_5}
	\norm{v_t}_{L^\infty([-1,1])} \leq 1 + \norm{A_{-1}^t + A_1^t}_{L^\infty([-1,1])}. 
	\end{align}
	If $x \in [-t/4,t/4]$, then we have the trivial bound
	\begin{align*}
	f(x/t)\geq \min_{[-1/4,1/4]}f\geq 2m.
	\end{align*}
	Therefore,
	\begin{align}\label{pf_lem_ext&cut_1}
	\max_{[-t/4,t/4]}\{A_{-1}^t + A_1^t\} \leq \frac{\norm{f}_{\C^0}}{2m}\quad \forall 1\leq t\leq 4.
	\end{align}
	
	Suppose now $x\in (t/4,1)$. Recalling that $f$ is non-increasing in $[1/4,1]$, we deduce that $f(x/t)\geq f(1/t)$. Therefore,
	\begin{align}\label{pf_lem_ext&cut_2}
	\max_{(t/4,1)}A_1^t \leq 1\quad \forall 1\leq t\leq 4.
	\end{align}
	Similarly, using $f(z/t)\leq f(x/t)$ for every $x\leq z\leq 1$, we obtain
	\begin{align*}
	1-\omega_t(x) = 1- \frac{\int_{-1}^x f(z/t)\,dz}{\int_{-1}^1f(y/t)\,dy} = \frac{\int_{x}^{1}f(z/t)\,dz}{\int_{-1}^1f(y/t)\,dy} \leq \frac{f(x/t)(1-x)}{\int_{-1}^1f(y/t)\,dy}\leq \frac{f(x/t)}{\int_{-1}^1f(y/t)\,dy},
	\end{align*}
	and, taking into account \eqref{eqqq_1}, this implies
	\begin{align}\label{pf_lem_ext&cut_3}
	\max_{(t/4,1)}A_{-1}^t \leq \frac{\norm{f}_{\C^0}}{m}\quad \forall 1\leq t\leq 4.
	\end{align}
	In a completely analogous way, one shows that 
	\begin{align}\label{pf_lem_ext&cut_4}
	\max_{(-1,-t/4)}A_{-1}^t \leq 1\quad \text{and}\quad \max_{(-1,-t/4)}A_1^t \leq \frac{\norm{f}_{\C^0}}{m}\quad \forall 1\leq t\leq 4.
	\end{align}
	Combining \eqref{pf_lem_ext&cut_5}, \eqref{pf_lem_ext&cut_1}, \eqref{pf_lem_ext&cut_2}, \eqref{pf_lem_ext&cut_3} and \eqref{pf_lem_ext&cut_4}, we infer that
	\begin{align*}
	\norm{v_t}_{L^\infty([-1,1])} \leq 2+\frac{\norm{f}_{\C^0}}{m}\leq \infty \quad \forall 1\leq t\leq 4.
\end{align*}	 
This proves that the $4$-dilation-cutoff $\mu_{(\cdot)}:[1,4]\to \pl_\infty^p(\R)$  is $\d$-Lipschitz and that the estimates \eqref{lem_dilCut_1} hold for almost every $1\leq t\leq 4$.

}
\end{proof}

\begin{lemma}[lipschitz continuity of the flattening-interpolation]\label{lem_flatt}
	Let $\mu=f\L\in \pl_\infty^\infty(\R)\subseteq \pl_\infty^p(\R)$ with $\spt\,f=[-1,1]$ and assume the existence of a constant $c>0$ such that $f\geq c$ almost everywhere in $[-1,1]$.  Then the $[-1,1]$-flattening-interpolation $\mu_{(\cdot)}=f_{(\cdot)}\L$ of $\mu$ is $\d$-Lipschitz and, in particular, 
	\begin{align*}
	|\dot f|^{L^p}(t)\leq \left(\frac{1}{2} + \norm{f}_{L^\infty}\right)2^\frac{1}{p}\quad\text{and}\quad |\dot\mu|^{W_\infty}(t)\leq  \max\left\{\frac{1}{c},4\right\}\left(\frac{1}{2}+\norm{f}_{L^\infty}\right)
\end{align*}	 
hold for almost every $0\leq t\leq 1$.
\end{lemma}

\begin{proof}
	A density of the $[-1,1]$-flattening-interpolation is $f_t = (1-t)f + t2^{-1}\chi_{[-1,1]}$.
	The function $t\mapsto f_t(x)$ is Lipschitz for almost every $x\in \R$ and its derivative is $\partial_t f_t = 2^{-1}\chi_{[-1,1]} - f$	for almost every  $0\leq t\leq 1$.
	Therefore $\partial_tf_t$ belongs to $L^\infty(\R)$ and it is supported on $[-1,1]$. Thus, 
	\begin{align*}
	\norm{\partial_tf_t}_{L^p} \leq \norm{\partial_t f_t}_{L^\infty}2^\frac{1}{p}\leq \left(\frac{1}{2} + \norm{f}_{L^\infty}\right)2^\frac{1}{p}.
	\end{align*}
This proves the first inequality of the statement.
	
	Consider the time-dependent Borel vector field $(t,x)\mapsto v_t(x)$ defined by 
	\begin{align*}
	v_t(x):= -\frac{\int_{-1}^x\partial_tf_t(z)\,dz}{f_t(x)}\chi_{[-1,1]}(x) \quad \forall (t,x)\in [0,1]\times \R.
	\end{align*}
	 Then, recalling that $f\geq c$ almost everywhere in $[-1,1]$, we have $f_t\geq \min\{c,2^{-1}\}$ almost everywhere in $[0,1]\times [-1,1]$ and therefore 
	\begin{align*}
	\norm{v_t}_{L^\infty(\mu_t)} \leq \frac{2\norm{\partial_tf_t}_{L^\infty}}{\min\{c,2^{-1}\}}\leq \max\left\{\frac{1}{c},4\right\}\left(\frac{1}{2}+\norm{f}_{L^\infty}\right).
	\end{align*}
	Since $\int_{-1}^1\partial_tf_t(z)\,dz=0$, arguing as in \eqref{eq_previous} one proves that $(\mu_{(\cdot)},v_{(\cdot)})$ is a weak solution of the continuity equation in $[0,1]$.	Whence the second inequality follows from Lemma \ref{lem_ineqA}.

\end{proof}

\textit{Proof of Theorem \ref{th_d1}.}
	Let $Q>0$ and $\mu=f\L\in\pl_\infty^p(\R)$ with $\spt\,\mu\subseteq [a^*,b^*]$ and $\norm{f}_{L^p}\leq Q$. Let $r:=(b^*-a^*)/2>1/3$ and consider $\alpha,\beta,M,a,b,\lambda$ as in Definition \ref{def_standardConnection} and $\rho^r$ as in \eqref{rho}. Observe that $\lambda > 6r > 2$ and $M=M(a^*,b^*)$ depends only on $a^*$ and $b^*$.
	
	Let $\mu_{(\cdot)}^1=f_{(\cdot)}^1\L:[0,1]\to \pl_\infty^p(\R)$ be the $\rho$-interpolution of $\mu$ and set $L=L(a^*,b^*):=\Lip(\rho^r)$. Define $\mu^1:=\mu_1^1$ and $f^1:=f\star \rho^r$. Then $\mu^1=f^1\L$. By virtue of Lemma \ref{lem_interp}, the curve $\mu^1_{(\cdot)}$ is $\d$-Lipschitz and satisfies the estimates 
	\begin{align}\label{eq_mu1}
	|\dot f^1|^{L^p}(t)\leq 2Q\quad\text{and}\quad |\dot \mu^1|^{W_\infty}(t)\leq 3r
	\end{align}	 
	for almost every $0\leq t\leq 1$, and the function $f^1$ enjoys the properties
	\begin{align}\label{eq_f1}
	\begin{matrix*}[l]
	f^1\in \C^\infty_c(\R), & \spt\,f^1=[a,b], & \Lip(f^1)\leq L\\
	f^1>0 \text{ in }(a,b), &f^1\uparrow \text{ in }(a,\alpha), & f^1\downarrow \text{ in }(\beta,b).
	\end{matrix*}
	\end{align}

	Let $\mu^2_{(\cdot)}=f^2_{(\cdot)}\L:[0,1]\to \pl_\infty^p(\R)$ be the $w$-shift of $\mu^1$ and set $\mu^2:=\mu_1^2$ and $f^2:=f^1(\cdot - w)$. Then $\mu^2 = f^2\L$. By \eqref{eq_f1} and Lemma \ref{lem_shift}, the curve $\mu^2_{(\cdot)}$ is $\d$-Lipschitz and satisfies the estimates
	\begin{align}\label{eq_mu2}
	|\dot f^2|^{L^p}(t)\leq LM(8r)^\frac{1}{p}\quad\text{and}\quad |\dot \mu^2|^{W_\infty}(t)\leq M
	\end{align}	 
	for almost every $0\leq t\leq 1$, and the function $f^2$ enjoys the properties
	\begin{align}\label{eq_f2}
	\begin{matrix*}[l]
	f^2\in \C^\infty_c(\R), & \spt\,f^2=\left[-\frac{\lambda}{2},\frac{\lambda}{2}\right], & \Lip(f^2)\leq L\\
	f^2>0 \text{ in }\left(-\frac{\lambda}{2},\frac{\lambda}{2}\right), &f^2\uparrow \text{ in }\left(-\frac{\lambda}{2},-\frac{\lambda}{2}+3r\right), & f^2\downarrow \text{ in }\left(\frac{\lambda}{2}-3r,\frac{\lambda}{2}\right).
	\end{matrix*}
	\end{align}

Let $\mu^3_{(\cdot)}=f^3_{(\cdot)}\L:[0,1]\to \pl_\infty^p(\R)$ the $(2/\lambda)$-scaling of $\mu^2$ and set $\mu^3:=\mu^3_1$ and $f^3:=\frac{\lambda}{2} f^2(\frac{\lambda}{2}\cdot)$. Then, by \eqref{eq_f2}, Lemma \ref{lem_scaling} and recalling that $\lambda \leq 8r$, the curve $\mu^3_{(\cdot)}$ satisfies the estimates
\begin{align}\label{eq_mu3}
	|\dot f^3|^{L^p}(t)\leq  32 Lr^2(4r-1)(8r)^\frac{1}{p} \quad \text{and}\quad |\dot \mu^3|^{W_\infty}(t) \leq 4r-1
\end{align}
for almost every $0\leq t\leq 1$ and the function $f^3$ enjoys the properties
\begin{align}
\label{eq_f3}
\begin{matrix*}[l]
f^3\in\C^\infty_c(\R),& \spt\,f^3= \left[-1,1\right],& \norm{f^3}_{\C^0}\leq \Lip(f^3)\leq \frac{\lambda^2 L}{4}, 
\\
f^3>0 \text{ in }\left(-1,1\right),&  f^3\uparrow \text{ in }\left(-1,-\frac{1}{4}\right), & f^3\downarrow \text{ in }\left(\frac{1}{4},1\right).
\end{matrix*}
\end{align}

Define now the functions $\tilde{f}\in L^p_c(\R)$ and $\tilde\rho^r\in \C^\infty_c(\R)$ as
\begin{align*}
\tilde{f}(x):= \frac{\lambda}{2}f\left(\frac{\lambda}{2}x+w\right),\quad \tilde{\rho}^r(x):= \frac{\lambda}{2}\rho^r\left(\frac{\lambda}{2}x\right).
\end{align*}
It is easy to verify that $f^3= \tilde{f}\star \tilde\rho^r$. Moreover, recalling the definition of $r$ and the assumption $[\alpha,\beta]\subseteq[a^*,b^*]$, it follows that
\begin{align*}
\spt \,\tilde f = \left[-1+\frac{6r}{\beta-\alpha + 6r},1-\frac{6r}{\beta-\alpha +6r}\right]\subseteq \left[-\frac{1}{4},\frac{1}{4}\right].
\end{align*}
 Therefore, by monotonicity and symmetry of $\tilde{\rho}^r$, and recalling that $\tilde f$ is a probability density, for every $x\in [-1/4,1/4]$ we have
\begin{align*}
f^3(x) \geq  \int_{-1/4}^{1/4}\tilde{f}(z)\tilde{\rho}^r(x-z)\,dz \geq \tilde{\rho}^r\left(\frac{1}{2}\right)= \frac{\lambda}{2}\rho^r\left(\frac{\lambda}{4}\right).
\end{align*}
Since $6r\leq \lambda \leq 8r$ and $\rho^r >0$ in $(-3r,3r)$, defining $m=m(a^*,b^*):= \frac{3}{2}r\rho^r(2r)$ we have
\begin{align}\label{eq_1}
\min_{[-1/4,1/4]}f^3 \geq 2m >0.
\end{align}

	Let $\mu^4_{(\cdot)}=f^4_{(\cdot)}\L:[1,4]\to \pl_\infty^p(\R)$ the $4$-dilation-cutoff of $	\mu^3$ and set $\mu^4:=\mu^4_4$ and $f^4:\R\to [0,\infty)$ be a density for $\mu^4$.
 Then, by \eqref{eq_f3}, \eqref{eq_1}, Lemma \ref{lem_dilCut} and using $\lambda \leq 8r$, the curve $\mu^4_{(\cdot)}$ is $\d$-Lipschitz and satisfies the estimates
\begin{align}\label{eq_mu4}
	|\dot f^4|^{L^p}(t)\leq 2^{10}\left(\frac{r^2 L}{m}\right)^2 2^\frac{1}{p}\quad \text{and}\quad |\dot \mu^4|^{W_\infty}(t)\leq 2+\frac{16r^2 L}{m}
\end{align}
for almost every $1\leq t\leq 4$ and the function $f^4$ satisfies the properties
\begin{align}
\begin{gathered}
\label{eq_f4}
 f^4> \frac{m}{2} \text{ almost everywhere in }\left[-1,1\right],\quad \norm{f^4}_{L^\infty}\leq \frac{16 r^2 L}{m}.
\end{gathered}
\end{align}

Let $\mu^5_{(\cdot)}=f^5_{(\cdot)}\L:[0,1]\to \pl_\infty^p(\R)$ the $[-1,1]$-flattening-interpolation of $\mu^4$. 
 Then, by \eqref{eq_f4} and Lemma \ref{lem_flatt}, the curve $\mu^5_{(\cdot)}$ is $\d$-Lipschitz and satisfies the estimates
\begin{align}\label{eq_mu5}
	|\dot f^5|^{L^p}(t)
	\leq \left(\frac{1}{2}+\frac{16r^2L}{m}\right)2^\frac{1}{p}
	\quad \text{and}\quad 
	|\dot \mu^5|^{W_\infty}(t)
	\leq \left(\frac{1}{2} + \frac{16 r^2 L}{m}\right)\max\left\{\frac{2}{m},4\right\}.
\end{align}
for almost every $0\leq t\leq 1$.

Finally, let $\mu_{(\cdot)}=f_{(\cdot)}\L:[0,7]\to \pl_\infty^p(\R)$ be the standard connection of $\mu$ with $[-1,1]$. Since $\mu$ is the concatenation of $\mu^1_{(\cdot)},\dots,\mu^5_{(\cdot)}$, then it immediately follows from \eqref{eq_mu1}, \eqref{eq_mu2}, \eqref{eq_mu3}, \eqref{eq_mu4} and \eqref{eq_mu5} that, with the current choice of $L,M$ and $m$ (which depend solely on $a^*$ and $b^*$),
\begin{align*}
|\dot f|^{L^p}(t)\leq C_{a^*,b^*,Q}^{L^p}(t) \quad \text{and} \quad |\dot \mu|^{W_\infty}(t)\leq C^{W_\infty}_{a^*,b^*}(t)
\end{align*}
hold
for almost every $0\leq t\leq 7$. Hence the theorem is proved.

{\flushright\qedsymbol

}

\section[proof of theorem {\nsiz 1.1}]{Proof of Theorem \ref{th_dn}}
\renewcommand{\L}{\mathscr{L}}

Before commencing the proof of Theorem \ref{th_dn}, let us introduce some notation. The $n$-dimensional Euclidean space $\Rn$ will be identified with the product $\R^{n-1}\times \R$, and its generic element will be denoted by $\Rn\ni\mathbf{x}=(x,y)\in  \R^{n-1}\times \R$. For any $\ell>0$, $Q_\ell:=[-\ell,\ell]^n$ denotes the $n$-dimensional cube centered in the origin and with side of length $2\ell$. 

For any function $f\in L^1(\Rn)$, we define the \emph{$\R^{n-1}$-projection of $f$}  and the \emph{$x$-residue of $f$}, denoted by $\pi(f):\R^{n-1}\to [0,\infty]$ and  $f^x:\R\to [0,\infty]$ respectively, by fixing a Borel-measurable representative of $f$ in the class $L^1(\Rn)$ and setting
\begin{gather*}
\pi(f)(x):=\int_\R f(x,z)\,dz,\quad 
f^x(y):=
\begin{cases}
\frac{f(x,y)}{\pi(f)(x)} &,\text{ if } \pi(f)(x)>0\\
\chi_{[0,1]}(y) &,\text{ if }\pi(f)(x)=0
\end{cases}.
\end{gather*}

Recall that if $\mu\in \P(\R^{k+m})$ is a probability measure and $F:\R^{k+m}\to \R^k$, $1\leq k <n$, is a Borel map, then the disintegration theorem for measures in Polish spaces (see e.g. \cite{Bogachev2007} or the appendix of \cite{FigalliGlaudo2021}) there exists a family $(\mu_F^\xi)_{\xi \in \R^{k}}\subseteq \P(\R^{k+m})$ such that:
\begin{enumerate}
\item $\xi\mapsto \mu_F^\xi$ is Borel-measurable with respect to the topology of the narrow convergence in $\P(\R^{k+m})$ (equivalently, $\xi\mapsto \mu_F^\xi(B)$ is Borel-measurable for every $B\in \B(\R^{k+m})$),
\item $\spt\,\mu_F^\xi\subseteq F^{-1}(\xi)$ for $F_\#\mu$-almost every $\xi\in \R^k$, and
\item  the measure $F_\#\mu\otimes (\mu_F^\xi)_\xi$ satisfies the property
\begin{align*}
\int_{\R^{k+m}} \psi(\BF x)\,d\mu(\BF x) = \int_{\R^k}\int_{F^{-1}(\xi)} \psi(\zeta)\,d\mu_F^\xi(\zeta)\,dF_\#\mu(\xi),
\end{align*}
for every non-negative Borel-measurable $\psi:\R^{k+m}\to[0,\infty]$.
\end{enumerate}  
The measure $F_\#\mu\otimes (\mu_F^\xi)_\xi$ is called \emph{Borel disintegration of $\mu$ with respect to $F$}.

\begin{lemma}\label{lem_existenceN_disintegration}
	Let $\mu=f\Ln\in \pl_\infty^p(\Rn)$, with $1\leq p\leq \infty$. For any choice of the representative $f$ in $L^1(\Rn)$ we have:
	\begin{enumerate}
	\item \label{lem_exN_dis_i} $(\pi(f)\L^\m)\otimes (f^x\L^1)_x$ is a Borel disintegration of $\mu$ with respect to the projection $(x,y)\mapsto x$, if we identify $\spt\,f^x\simeq \{x\}\times \spt \,f^x$;
	\item \label{lem_exN_dis_iii}$\norm{x\mapsto \pi(f)(x)\norm{f^x}_{L^p(\R)}}_{L^p(\R^\m)}= \norm{f}_{L^p(\Rn)}$;
	\item \label{lem_exN_dis_ii} $\pi(\mu):=\pi(f)\L^\m\in \pl_\infty^p(\R^\m)$ and $\mu^x:=f^x\L^1\in \pl_\infty^p(\R)$ for $\pi(\mu)$-almost every $x\in \R^\m$. In particular, if $\spt\,\mu\subseteq Q_r$, then
	\begin{align}\label{eq_recalling}
	\norm{\pi(f)}_{L^p(\R^\m)} \leq (2r)^\frac{p-1}{p}\norm{f}_{L^p(\Rn)},
	\end{align}
	where $(p-1)/p:=1$, if $p=\infty$.
\end{enumerate}		
\end{lemma}

\begin{proof}
	Let $\mu=f\Ln\in \pl_\infty^p(\Rn)$ and fix a Borel-measurable representative $f$ in the equivalence class in $L^1(\Rn)$. The functions $(x,y)\mapsto \pi(f)(x)$ and $(x,y)\mapsto f^x(y)$ are Borel-measurable. Therefore also $x\mapsto f^x\L^1$ is a Borel-measurable map $\R^\m\to \P(\R)$. For any function $g\in \C^0_b(\Rn)$ we have 
	\begin{align*}
	\int_\Rn g(\BF x)\,d\mu(\BF x) = \int_{\R^\m}\int_\R g(x,y)f(x,y)\,dy\,dx = \int_{\R^\m}\int_\R g(x,y)\,d\mu^x(y)\,d\pi(\mu)(x),
\end{align*}		
and clearly $\pi(\mu)\in \P(\R^\m)$ and $\mu^x\in \P(\R)$ for $\pi(\mu)$-almost every $x\in \R^\m$.  Therefore, $(\pi(f)\L^\m)\otimes (f^x\L^1)_x$ is a Borel disintegration of $\mu$ with respect to the projection $(x,y)\mapsto x$, once the support of $f^x$ is identified with $\{x\}\times \spt\,f^x$.
	
	By compactness of the support of $\mu$, there exists $r>0$ such that $\spt\,\mu = \spt\,f\subseteq Q_r$, and therefore $\spt \,\pi(\mu) \subseteq [-r,r]^\m$ and $\spt\,\mu^x\subseteq [-r,r]$ for $\pi(\mu)$-almost every $x\in \R^\m$. We have
	\begin{align}\label{eq_dis_1}
	\begin{split}
	\norm{f}_{L^p(\Rn)} 
	&= \norm{x\mapsto \pi(f)(x)\norm{f^x}_{L^p([-r,r])}}_{L^p([-r,r]^\m)}.
	\end{split}
	\end{align}
	This identity is immediate for $1\leq p<\infty$, while it is slightly less trivial for $p=\infty$ (see e.g. Lemma A.1 in \cite{AldrigoPLpq}). Since $f\in L^p(\Rn)$, $x\mapsto \pi(f)(x)\norm{f^x}_{L^p(\R)}$ is $p$-integrable in $\R^\m$, and therefore $\norm{f^x}_{L^p(\R)}<\infty$ for $\pi(\mu)$-almost every $x\in \R^\m$. This proves that $\mu^x\in \pl_\infty^p(\R)$ for $\pi(\mu)$-almost every $x\in \R^\m$.
	
	We now only need to show that $\pi(f)\in L^p(\R^\m)$. Suppose first $1\leq p <\infty$. By Jensen's inequality,
	\begin{align*}
	\begin{split}
	\norm{\pi(f)}_{L^p(\R^\m)} = \left(\int_{[-r,r]^\m}\left(\int_{-r}^r f(x,y)\,dy\right)^pdx\right)^\frac{1}{p}
	\leq (2r)^\frac{p-1}{p}\norm{f}_{L^p(\Rn)}.
	\end{split}
	\end{align*}
	Finally, for $p=\infty$, the first identity of \eqref{eq_dis_1} proves that $\norm{f(x,\cdot)}_{L^\infty(\R)}\leq \norm{f}_{L^\infty(\Rn)}$ for almost every $x\in \R^\m$. Therefore, for any of such $x$'s, we have
	\begin{align*}
	\pi(f)(x) = \int_{-r}^r f(x,y)\,dy \leq 2r \norm{f}_{L^\infty(\Rn)}.
	\end{align*}
	This yields $\norm{\pi(f)}_{L^\infty(\R^\m)}\leq 2r \norm{f}_{L^\infty(\Rn)}$.

\end{proof}

From now on we implicitly assume that whenever the notation $f\Ln = (\pi(f)\L^\m)\otimes(f^x\L^1)_x$ is used, a Borel-measurable representative of $f$ in $L^1(\Rn)$ is fixed such that the map $x\mapsto \mu^x$ is Borel-measurable.

\begin{lemma}\label{lem_existenceN_measurable_selection}
	Let $\mu_{(\cdot)}:[0,T]\to \P_\infty(\R^{k+m})$ be a narrowly continuous curve, $L>0$ and let $P:\R^{k+m}\to \R^k$ be the orthogonal projection onto the first $k$ coordinates. Suppose the existence of $\kappa\in \P_\infty(\R^k)$ and a Borel disintegration with respect to the projection $P$ of $\mu_t$ of the form $\kappa\otimes(\mu^\xi_t)_{\xi}$, with $\mu^\xi_{(\cdot)}:[0,T]\to \P_\infty(\R^{m})$ $W_\infty$-absolutely continuous and $|\dot \mu^\xi|^{W_\infty}(t)\leq L$ for almost every $0\leq t\leq T$ and $\kappa$-almost every $\xi\in \R^k$. Then $\mu_{(\cdot)}$ is $W_\infty$-absolutely continuous, and $|\dot\mu|^{W_\infty}(t)\leq L$ for almost every $0\leq t\leq T$.
\end{lemma}

\begin{proof}
	Fix $0\leq t\leq s\leq T$ and let $A:=\{(y,z)\in \R^m\times\R^m:|y-z|\leq L(s-t)\}$. Under the standing assumptions, $W_\infty(\mu_s^\xi,\mu_t^\xi)\leq L(s-t)$ for $\kappa$-almost every $\xi\in \R^k$. Consider the families
	\begin{align*}
	\Omega_\xi:=\left\{\omega\in \mathcal{P}(\R^m\times \R^m): (p_1)_\#\omega = \mu_s^\xi,\,(p_2)_\#\omega =\mu_t^\xi,\,\omega(A)=1\right\}\subseteq \mathcal{M}(\R^m\times \R^m),
	\end{align*}
	where $p_1,p_2:\R^m\times \R^m\to \R^m$ are the projections $p_1(y,z):=y$ and $p_2(y,z):=z$. Observe that $\Omega_\xi$ is a subset of $\Gamma(\mu_s^\xi,\mu^\xi_t)$ for $\kappa$-almost every $\xi\in\R^k$. Because $\Omega_\xi$ always contains the set of $\infty$-optimal couplings $\Gamma_\infty(\mu_s^\xi,\mu_t^\xi)\neq \varnothing$ (see e.g. \cite{champion2008wasserstein,givens}), it is non-empty for $\kappa$-almost every $\xi\in \R^k$. We claim that the set
	\begin{align*}
	\Graph(\Omega):=\left\{(\xi,\omega)\in \R^k\times\mathcal{M}(\R^m\times\R^m): \omega\in \Omega_\xi\right\}
	\end{align*}
	is Borel, if $\mathcal{M}(\R^m\times \R^m)$ is endowed with the topology of the narrow convergence. First, observe that $(p_1)_\#,(p_2)_\#:\mathcal{M}(\R^m\times\R^m)\to \mathcal{M}(\R^m)$ are narrowly continuous maps, and since $A\subseteq \R^m\times \R^m$ is closed, then $\omega\mapsto \omega(A)$ is narrowly upper semicontinuous.	
	 Therefore, the family $\Omega_\xi$ is a closed subset of $\mathcal{P}(\R^m\times \R^m)$ for $\kappa$-almost every $\xi\in \R^k$. Recall that maps $\xi\mapsto \mu^\xi_s$ and $\xi\mapsto \mu^\xi_t$ are both Borel. This implies that the map $\Phi:\R^k\times \mathcal{M}(\R^m\times \R^m)\to \mathcal{M}(\R^m)\times \mathcal{M}(\R^m)\times \R$ defined by
	 \begin{align*}
	 \Phi(\xi,\omega):= ((p_1)_\#\omega - \mu^\xi_s,(p_2)_\# \omega - \mu^\xi_t, \omega(A) -1 ) 
	 \end{align*}
	is Borel. Since $\Graph(\Omega)=\Phi^{-1}(\{0\})$, then $\Graph(\Omega)$ is Borel.
	
	Thanks to classical measurable selection theorems (see e.g. Jankov--von Neumann's theorem), there exists a $\kappa$-measurable map $\xi\mapsto \omega_\xi$ with $\omega_\xi\in \Omega_\xi$ for $\kappa$-almost every $\xi\in \R^k$. Define now the probability measure $\gamma \in \P(\R^{k+m}\times \R^{k+m})$ by integration of continuous and bounded functions $\psi\in \C^0_b(\R^{k+m}\times\R^{k+m})$ as
	\begin{align*}
	\int_{\R^{k+m}\times \R^{k+m}} \psi((\xi,y),(\zeta,z))\,d\gamma((\xi,y),(\zeta,z))
	:= \int_{\R^k}\int_{\R^m\times \R^m}\psi((\xi,y),(\xi,z))\,d\omega_\xi(y,z)\,d\kappa(\xi).
\end{align*}		
	Then $\gamma\in \Gamma(\mu_s,\mu_t)$ and $\gamma(\{(\BF x, \BF y)\in \R^{k+m}\times \R^{k+m}:|\BF x - \BF y|\leq L(s-t)\})=1$. Therefore $W_\infty(\mu_s,\mu_t)\leq L(s-t)$. The claim follows by arbitrariness of the choices of $0\leq t\leq s\leq T$.

	\end{proof}
	
\begin{lemma}\label{lem_BorelMeas}
	Let $\mu \in \pl_\infty^p(\Rn)$ with $1\leq p\leq \infty$, $\spt\,\mu\subseteq Q_r$, and let $\pi(\mu)\otimes (\mu^x)_x$ be a Borel disintegration of $\mu$. Denote by $\mu^x_{(\cdot)}:[0,7]\to \pl_\infty^p(\R)$ be the $(-r,r,\rho^r)$-standard connection of $\mu^x$ with $[-1,1]$. Then $\mu_t:=\pi(\mu)\otimes (\mu^x_t)_x$ is a Borel measure for every $0\leq t\leq 7$.
\end{lemma}

\begin{proof}
	For each $\mu^x$, let $\alpha^x,\beta^x,r^x,M^x,a^x,b^x,w^x,\lambda^x$ be the constants $\alpha,\beta,r,M,a,b,$ $w,\lambda$ appearing in Definition \ref{def_standardConnection} relative to $\mu^x$, and let $\mu^{1,x}_{(\cdot)},\mu^{2,x}_{(\cdot)}\mu^{3,x}_{(\cdot)}\mu^{4,x}_{(\cdot)}\mu^{5,x}_{(\cdot)}$ be respectively the $\rho$-interpolution of $\mu^x$, the $w^x$-shift of $\mu^{1,x}_1$, the $(2/\lambda^x)$-scaling of $\mu^{2,x}_1$, the $4$-dilation-cutoff of $\mu^{3,x}_1$ and the $[-1,1]$-flattening-interpolation of $\mu^{4,x}_4$. We show that 
	$\mu^{j}_t:=\pi(\mu) \otimes (\mu^{j,x}_t)_x$  is  Borel for every $1 \leq j\leq 5$ and every $t$ for which $\mu^{j,x}_t$ is defined.
	
	For $j=1$, observe that 
	\begin{align*}
	\mu_t^1 = \pi(\mu)\otimes (\mu^x\star \kappa_t),\quad \kappa_t := (1-t)\delta_0 + \rho^r(\cdot/t)\L^1\quad \forall 0\leq t\leq 1.
\end{align*}	 
	Since $x\mapsto \mu^x(B)$ is Borel-measurable for every $B\in \B(\R)$, then so it is $x\mapsto \mu^x\star \kappa_t (B) = (1-t)\mu^x(B) +  \int_B f^x\star\rho^r(\cdot/t))\,d\L^1$, and this is enough to guarantee the Borel-measurability of $\mu^1_t$.
	
	A sufficient condition for the Borel-measurability of $\mu^j_t$, for $j=2$ and $j=3$, is the Borel-measurability of $x\mapsto w^x$ and $x\mapsto \lambda^x$ respectively. Since $w^x= (a^x+b^x)/2$ and $\lambda^x = b^x-a^x$, then it will be enough to show that $x\mapsto a^x$ and $x\mapsto b^x$ are both Borel-measurable. Let $f^{1,x}=f^x\star \rho^r$ be the density of $\mu^{1,x}_1$. Then $(x,y)\mapsto f^{1,x}(y)$ is Borel-measurable and $y\mapsto f^{1,x}(y)$ is smooth for $\pi(\mu)$-almost every $x\in \R^\m$. Observe that $[a^x,b^x] = \spt\,f^{1,x}$, i.e. 
	$a^x = \inf\,\spt\,f^{1,x}$ and  $b^x = \sup\,\spt\,f^{1,x}$.
Therefore, for any constant $c\in \R$, we have that
\begin{align*}
	\left\{x\in \R^\m: a^x <c\right\} 
	= \bigcup_{\substack{q\in \BF Q\\ q<c}} \left\{x\in \R^\m : f^{1,x}(q)>0\right\},
\end{align*}
is the countable union of Borel subsets; hence it is Borel. Similarly, 
\begin{align*}
\left\{x\in \R^\m : b^x>c\right\} = \bigcup_{\substack{q\in \BF Q\\ q>c}} \left\{x\in \R^\m: f^{1,x}(q)>0\right\}.
\end{align*}
Thus $x\mapsto a^x$ and $x\mapsto b^x$ are both Borel-measurable.

Finally, for $j=4$ and $j=5$ the proof is trivial.

\end{proof}

\textit{Proof of Theorem \ref{th_dn}.} Arguing as in Corollary \ref{cor}, it is enough to prove the following claim.

\begin{equation}\tag{$\claim$}\label{claim}
\begin{minipage}{0.85\textwidth}
\itshape “For any $\mu\in \pl_\infty^p(\Rn)$, with $1\leq p\leq \infty$ and $n\geq 1$ there exists a $\d$-absolutely continuous curve $\mu_{(\cdot)}:[0,1]\to \pl_\infty^p(\Rn)$ such that $\mu_0=\mu$ and $\mu_1=2^{-n}\Ln\llcorner Q_1$, where $Q_1:=[-1,1]^n$.” 
\end{minipage}
\end{equation}

We prove \eqref{claim} by induction on the dimension. If $n=1$, then \eqref{claim} follows from Theorem \ref{th_d1}. Suppose now $n\geq 2$ and the claim to hold true up to dimension $n-1$. Write 
\begin{align*}
\mu= f\Ln = (\pi(f)\L^\m)\otimes (f^x\L^1)_x = \pi(\mu)\otimes (\mu^x)_x,
\end{align*}
where $f\in L^p(\Rn)$ is chosen in such a way that $x\mapsto \mu^x$ is Borel-measurable. Then, by Lemma \ref{lem_existenceN_disintegration}, $\pi(\mu)\in \pl_\infty^p(\R^\m)$ and $\mu^x\in \pl_\infty^p(\R)$ for $\pi(\mu)$-almost every $x\in \R^\m$. Let $Q_r\subseteq \Rn$ be such that $\spt\,\mu\subseteq Q_r$, so that $\spt \,\mu^x\subseteq [-r,r]$ for $\pi(\mu)$-almost every $x$. For any of such $x$, let $\mu^x_{(\cdot)}=f^x_{(\cdot)}\L^1:[0,7]\to \pl_\infty^p(\R)$ be the $(-r,r,\rho^r)$-standard connection with $[-1,1]$. Without loss of generality, we can assume $r>1/3$. Then, by Theorem \ref{th_d1}, it follows that $\mu^x_{(\cdot)}$ is $\d$-Lipschitz, and in particular we have the estimates
\begin{gather}\label{eqq_1}
|\dot f^x|^{L^p}(t)\leq C^{L^p}_{x}(t):=
\begin{cases}
	 2\norm{f^x}_{L^p(\R)}																&, \text{ if } 0\leq t < 1\\
	 LM(8r)^\frac{1}{p}												&, \text{ if } 1\leq t < 2\\
	 32 Lr^2(4r-1)(8r)^\frac{1}{p}									&, \text{ if } 2\leq t < 3\\
	 2^{10}\left(\frac{r^2 L}{m}\right)^2 2^\frac{1}{p}				 		&, \text{ if } 3\leq t < 6\\
	 \left(\frac{1}{2}+\frac{16r^2L}{m}\right)2^\frac{1}{p}				&, \text{ if } 6\leq t \leq 7
\end{cases},
\end{gather}
\begin{gather}
\label{eqq_2}
|\dot \mu^x|^{W_\infty}(t)\leq C^{W_\infty}(t):=
\begin{cases}
	 3r 																		&, \text{ if } 0\leq t < 1\\
	 M																		&, \text{ if } 1\leq t < 2\\
	4r-1																		&, \text{ if } 2\leq t < 3\\
	2+\frac{16 r^2 L}{m} 														&, \text{ if } 3\leq t < 6\\
	\left(\frac{1}{2} + \frac{16 r^2 L}{m}\right)\max\left\{\frac{2}{m},4\right\}	&, \text{ if } 6\leq t \leq 7,
\end{cases}
\end{gather}
where $L,M,m$ depend only on $r$.

Define $\mu_t:=\pi(\mu)\otimes (\mu^x_t)_x$. Then Lemma \ref{lem_BorelMeas} ensures that $\mu_t$ is a Borel measure for every $0\leq t\leq 7$, and (a Borel-measurable choice of) the density of $\mu_t$ is 
\begin{align*}
f_t(x,y):= \pi(f)(x)f^x_t(y).
\end{align*}
If $1\leq p <\infty$, then, by Minkowski and Lemma \ref{lem_existenceN_disintegration}
\begin{align*}
\begin{split}
\norm{f_s-f_t}_{L^p} 
	&\leq \left(\int_{\R^\m} \left(\int_t^s \pi(f)(x)C_x^{L^p}(\tau)\,d\tau\right)^p \,dx\right)^\frac{1}{p}\\
	&\leq \int_t^s \norm{x\mapsto \pi(f)(x)C_x^{L^p}(\tau)}_{L^p(\R^\m)}\,d\tau,
\end{split}
\end{align*}
and by Lemma \ref{lem_existenceN_disintegration}, the function $x\mapsto \pi(f)(x)C_x^{L^p}(\tau)$ belongs to $L^p(\R^\m)$ and, recalling \eqref{eq_recalling} and \eqref{eqq_1} we obtain
\begin{align*}
\norm{x\mapsto \pi(f)(x)C_x^{L^p}(\tau)}_{L^p(\R^\m)} \leq
\begin{cases}
	 2\norm{f}_{L^p}																&, \text{ if } 0\leq \tau < 1\\
	 2 LMr 4^\frac{1}{p}	\norm{f}_{L^p}													&, \text{ if } 1\leq \tau < 2\\
	 32 Lr^2(4r-1)(8r)^\frac{1}{p}(2r)^\frac{p-1}{p}	\norm{f}_{L^p}										&, \text{ if } 2\leq \tau < 3\\
	 2^{10}\left(\frac{r^2 L}{m}\right)^2 2^\frac{1}{p}	(2r)^\frac{p-1}{p}	\norm{f}_{L^p}				 		&, \text{ if } 3\leq \tau < 6\\
	 \left(\frac{1}{2}+\frac{16r^2L}{m}\right)2^\frac{1}{p}	\norm{f}_{L^p}				&, \text{ if } 6\leq \tau \leq 7
\end{cases}.
\end{align*}
Thus, $f_{(\cdot)}:[0,7]\to L^p(\Rn)$ is $L^p$-absolutely continuous and has bounded $L^p$-metric derivative.

If $p=\infty$, then a similar argument
proves that there exists a constant $\tilde C$ such that 
$\norm{x\mapsto \pi(f)C_x^{L^\infty}(\tau)}_{L^\infty}\leq \tilde{C}$ for every $0\leq \tau\leq 7$. Whence it follows that $f_{(\cdot)}:[0,7]\to \pl_\infty^\infty(\Rn)$ is $L^\infty$-absolutely continuous and has bounded $L^\infty$-metric derivative.

Since the function $C^{W_\infty}$ of \eqref{eqq_2} does not depend on $x$, and $L^p$-absolute continuity of a curve of equi-compactly supported densities implies narrow continuity, Lemma \ref{lem_existenceN_measurable_selection} applies and $\mu_{(\cdot)}:[0,7]\to \pl_\infty^p(\Rn)$ is $W_\infty$-Lipschitz with $W_\infty$-metric derivative dominated by $C^{W_\infty}$.  Therefore, $\mu_{(\cdot)}$ is $\d$-Lipschitz. 

Define $\tilde{\mu}:=\mu_7 = \pi(\mu)\otimes(\frac{1}{2}\L^1\llcorner[-1,1])$. By induction, there exists a $\d$-Lipschitz curve $\gamma_{(\cdot)}=g_{(\cdot)}\L^\m:[0,1]\to \pl_\infty^p(\R^\m)$ such that $\gamma_0 = \pi(\mu)$ and $\gamma_1=2^{-(n-1)}\L^\m\llcorner [-1,1]^\m$. Let $\nu_{(\cdot)}:[0,1]\to \pl_\infty^p(\Rn)$ be the curve
\begin{align*}
\nu_t = \gamma_t\otimes \left(\frac{1}{2}\L^1\llcorner [-1,1]\right) = G_t\Ln,\quad G_t(x,y)=\frac{g_t(x)\chi_{[-1,1]}(y)}{2}.
\end{align*}
Then, by Lemma \ref{lem_existenceN_disintegration}, 
\begin{align*}
\norm{G_s-G_t}_{L^p} = 2^\frac{1}{p}\frac{\norm{g_s-g_t}_{L^p(\R^\m)}}{2},
\end{align*}
with the convention $c^\frac{1}{\infty}:=1$ for every $c>0$. Thus $G_{(\cdot)}$ inherits from $g_{(\cdot)}$ a bounded $L^p$-metric derivative. 

On the other hand, since the conditional measures of $\nu$ with respect to the projection onto the first $n-1$ coordinates is constantly equal to $\frac{1}{2}\L^1\llcorner[-1,1]$, Lemma \ref{lem_existenceN_measurable_selection}, combined with the $W_\infty$-Lipschitz continuity of $\gamma_{(\cdot)}$, guarantee the $W_\infty$-Lipschitz continuity of $\nu_{(\cdot)}$. Concatenating $\mu_{(\cdot)}$ with $\nu_{(\cdot)}$ and reparametrizing in time, we obtain a $\d$-absolutely continuous curve that connects $\mu$ with $2^{-n}\Ln\llcorner Q_1$. Hence, by induction, \eqref{claim} holds true.

\hfill \qedsymbol

\section{Final remarks}

Let us recall that in the classical Wasserstein spaces $(\P_q(\Rn),W_q)$, with $1<q\leq\infty$, Benamou--Brenier's theorem (see e.g. \cite{AldrigoPLpq,ambrosio2021,BenamouBrenier2000}) implies that the metric $W_q$ coincides with the \emph{intrinsic length distance $\hat W_q$ associated with $W_q$}, given by
\begin{align}\notag
\hat W_q(\mu,\nu)&:=\min\left\{\int_0^1 |\dot \pi|^{W_q}(t)\,dt:\,\begin{matrix}
\pi_{(\cdot)}\in \ac([0,1];\P_q(\Rn))\hfill\\
\pi_0=\mu\hfill\\
\pi_1=\nu\hfill
\end{matrix}
\right\}\\ 
&=\min\left\{\int_0^1 \norm{v_t}_{L^q(\pi_t)}\,dt:\,
\begin{matrix}
\partial_t \pi_t + \div(v_t\pi_t) = 0\hfill\\
\pi_0=\mu\hfill\\
\pi_1=\nu\hfill
\end{matrix}\notag 
\right\}.
\end{align}
The equality $W_q=\widehat W_q$ is a key reason why gradient flows and minimizing movements are naturally formulated in Wasserstein spaces, as it allows one to work directly with absolutely continuous curves satisfying the continuity equation, rather than only with optimal plans; see \cite{ambrosio2005gradient,JKO,Otto2001,Santambrogio2017}. It is also closely related to the role of displacement/geodesic convexity in optimal transport, starting from \cite{McCann1997}, as explained in Chapter 11 of \cite{ambrosio2005gradient}.

In this paper we have proved that for every choice of $\mu,\nu\in \pl_\infty^p(\Rn)$, the family of $\d$-absolutely continuous curves $\gamma_{(\cdot)}:[0,1]\to \pl_\infty^p(\Rn)$ with $\gamma_0=\mu$ and $\gamma_1=\nu$ is non-empty. This naturally leads to the definition of the {intrinsic length distance associated with $\d$} $\D:\pl_\infty^p(\Rn)\times \pl_\infty^p(\Rn)\to [0,\infty)$
\begin{align*}
\D(\mu,\nu):=\inf\left\{\int_0^1|\dot\gamma|(t)\,dt:\,\begin{matrix} \gamma_{(\cdot)}\in \ac([0,1];\pl_\infty^p(\Rn))\hfill\\
\gamma_0=\mu\hfill\\
\gamma_1=\nu\hfill\end{matrix}\right\}.
\end{align*}
Taking into account Theorem 5.2 and Theorem 5.3 of \cite{AldrigoPLpq}, we can also write the Eulerian formulations
\begin{gather*}
\D(\mu,\nu)=\inf\left\{\int_0^1\left(\norm{v_t}_{L^\infty(\gamma_t)} + \norm{\div(v_th_t)}_{L^p}\right)\,dt:\,
\begin{matrix}
\gamma_t = h_t\Ln\hfill \\
\partial_th_t + \div(v_th_t) = 0 \hfill\\
h_0\Ln = \mu \hfill\\
h_1\Ln = \nu \hfill
\end{matrix}
\right\}
\end{gather*}
if $1< p\leq \infty$, and
\begin{gather*}
\hat{\mathfrak{d}}_\infty^1(\mu,\nu) = \inf\left\{\int_0^1\left(\norm{v_t}_{L^\infty(\gamma_t)} + \norm{\div(v_th_t)}_{\tv}\right)\,dt:\,
\begin{matrix}
\gamma_t = h_t\Ln\hfill \\
\partial_th_t + \div(v_th_t) = 0 \hfill\\
h_0\Ln = \mu \hfill\\
h_1\Ln = \nu \hfill
\end{matrix}
\right\}
\end{gather*}
where $\norm{\div(v_th_t)}_{L^p}$ (resp. $\norm{\div(v_th_t)}_{\tv}$, if $p=1$) is set to be infinity, if the function $(t,x)\mapsto v_t(x)f_t(x)$ does not admit a weak $x$-divergence in $L^1_tL^p_x([0,1]\times \Rn)$ (resp. in $\mathcal{M}^1_t\text{FV}_x([0,1]\times \Rn)$). 

It is elementary to show that $\D$ is in fact a metric on $\pl_\infty^p(\Rn)$, $(\pl_\infty^p(\Rn),\D)$ is a \emph{length space} and that the inequality $\d(\mu,\nu)\leq \D(\mu,\nu)$ holds for every choice of $\mu,\nu\in \pl_\infty^p(\Rn)$. Moreover, there are examples for which the inequality is strict. Therefore it is not possible to obtain a Benamou--Brenier-like theorem in the space $\pl_\infty^p(\Rn)$ (take, for instance, probabilities $\mu=f\L^1,\nu=g\L^1\in \pl_\infty^p(\R)$ with densities $f,g\in \Lip_c(\R)$, where $g=f(\cdot - w)$ for some large $|w|$).

On the other hand, two questions remain open and are left for future research.
\begin{enumerate}[label = \arabic*.]
\item Do $\d$ and $\D$ induce equivalent metric structures on $\pl_\infty^p(\Rn)$?
\item Is $(\pl_\infty^p(\Rn),\D)$ a \emph{geodesic metric space}?
\end{enumerate}

A positive answer to the first question would imply that the theory developed in Section 3 of \cite{AldrigoPLpq} for the existence of isoperimetric minimizing movements in $(\pl_\infty^\infty(\Rn),\d)$ remains valid in $(\pl_\infty^\infty(\Rn),\D)$. If the second question also has a positive answer, then the main tools of Otto calculus and of the theory of geodesically convex functionals become available in the $\pl_\infty^p$ setting as well. This could lead to a broader study of $\pl$-gradient flows for functionals depending on ratios of integral norms of a density and its derivatives.

\color{black}

\end{document}